\documentclass[11pt,a4paper]{article}

\usepackage{amsmath}
\usepackage{amsfonts}
\usepackage{amsthm}
\usepackage{amssymb}
\usepackage{url}
\usepackage[usenames,dvips]{color}
\usepackage{stmaryrd}
\usepackage[dvips]{graphicx}
\usepackage{psfrag, float}  

\DeclareMathAlphabet{\mathpzc}{OT1}{pzc}{m}{it}

\def\figlabel#1{\label{#1}}
\newtheorem{theorem}{Theorem}[section]

\newtheorem{corollary}[theorem]{Corollary}

\newtheorem{remark}[theorem]{Remark}

\newtheorem{conjecture}[theorem]{Conjecture}

\newcommand{\be}{\begin{equation}}
\newcommand{\ee}{\end{equation}}
\newcommand{\eps}{\varepsilon}
\newcommand{\ga}{\gamma}
\newcommand{\dps}{\displaystyle}
\newcommand{\RR}{\mathbb{R}}
\newcommand{\NN}{\mathbb{N}}
\newcommand{\CC}{\mathbb{C}}
\newcommand{\TT}{\mathbb{T}}
\newcommand{\ZZ}{\mathbb{Z}}

\newcommand{\LL}{\mathcal{L}}
\newcommand{\QQ}{\mathbb{Q}}
\newcommand{\YY}{\mathcal{Y}}

\newcommand{\DD}{\mathcal{D}}
\newcommand{\XX}{\mathcal{X}}
\newcommand{\OO}{\mathcal{O}}

\newcommand{\HH}{\mathcal{H}}
\newcommand{\VV}{\mathcal{V}}
\newcommand{\RRR}{\mathcal{R}}

\newcommand{\UU}{\mathcal{U}}

\newcommand{\AAA}{\mathcal{A}}

\newcommand{\CCC}{\mathcal{C}}

\newcommand{\xp}{x_p}
\newcommand{\yp}{y_p}

\newcommand{\ii}{^{-1}}
\newcommand{\de}{\delta}
\newcommand{\pa}{\partial}
\newcommand{\la}{\lambda}

\newcommand{\inn}{\mathrm{in}}
\newcommand{\out}{\mathrm{out}}

\newcommand{\kk}{\kappa}
\newcommand{\rr}{\rho}
\newcommand{\tro}{I}

\newcommand{\C}{c}

\newcommand{\tet}{\theta}
\newcommand{\ol}{\overline}

\renewcommand{\Re}{\mathrm{Re\, }}
\renewcommand{\Im}{\mathrm{Im\,}}
\newcommand{\wt}{\widetilde}
\newcommand{\wh}{\widehat}

\newcommand{\hmu}{\mu}

\begin{document}

\title{Splitting of separatrices in the resonances of nearly integrable Hamiltonian Systems of one and a half degrees of freedom}
\author{Marcel Guardia\thanks{\tt guardia@ias.edu, marcel.guardia@upc.edu}}
\maketitle
\medskip
\begin{center}$^\ast$
School of Mathematics\\
Institute for Advanced Study\\
Einstein Drive, Simonyi Hall\\
Princeton, New Jersey, 08540\end{center}

\begin{abstract}
In this paper we consider general nearly integrable analytic Hamiltonian systems of one and a half degrees of freedom which are a trigonometric polynomial in the angular state variable. In the resonances of these systems generically appear hyperbolic periodic orbits. We study the possible transversal intersections of their invariant manifolds, which is exponentially small, and we give an asymptotic formula for the measure of the splitting. We see that its asymptotic first order is of the form $K \varepsilon^{\beta} \text{e}^{-a/\varepsilon}$ and we identify the constants $K,\beta,a$ in terms of the system features. We compare our results with the classical Melnikov Theory and we show that, tipically, in the resonances of nearly integrable systems Melnikov Theory fails to predict correctly the constants $K$ and $\beta$ involved in the formula.

\end{abstract}

\section{Introduction}\label{sec:intro}
In this paper we consider nearly integrable analytic Hamiltonian systems of one and a half degrees of freedom. These Hamiltonians are of the form
\begin{equation}\label{def:HamResonant}
h(x,I,\tau;\de)=h_0(I)+\de h_1(x,I,\tau;\de),
\end{equation}
where $\de\ll 1$ is a small parameter, $(x,\tau)\in\TT^2$,
$I\in U\subset\RR$ and $h_0$ and  $h_1$ are analytic functions.

When $\de=0$, the Hamiltonian system is completely integrable (in
the sense of Liouville-Arnold)  and $I$ is a first integral of the system. That is, for $\de=0$ the phase space is foliated by
invariant tori $I=I^\ast$ with frequency vector $\omega(I^*)=(\pa_I h_0(I^*),1)$. The dynamics in these tori is quasiperiodic or periodic depending whether $\pa_I h_0(I^*)$ is rational or irrational. The natural question concerning nearly integrable Hamiltonian systems, is how the dynamics changes when one adds the perturbation. Namely, which tori persist and which break down, and, what new kinds of dynamics appear in the regions where the unperturbed invariant tori  break down.

The persistence of most of the invariant tori (in the measure sense) was shown by Kolmogorov-Arnol'd-Moser Theory (see \cite{Kolmogorov54,Arnold63a,Moser62}). The purpose of this paper is to study some particular dynamics in the complementary of this region. Namely, to study the dynamics in the region in which the tori that existed for the unperturbed system break down. These regions are usually called resonances and they appear when the frequency vector $\omega(I^*)$ of the unperturbed system is rationally dependent. This happens for $I=I^\ast$ such that $\pa_I h_0(I^*)\in \QQ$. We  call $n/m:1$-resonance to the region of the phase space $I\sim I^\ast$ such that $\pa_Ih_0(I^*)=n/m$, $n\in\ZZ$ and $m\in\NN$.  Under certain non-degeneracy conditions, in the resonances appear, among other invariant objects, low dimensional hyperbolic invariant tori which, for systems of one and a half degrees of freedom, are just hyperbolic periodic orbits. These objects have stable and unstable invariant manifolds, often called whiskers. The main goal of this paper is to study the possible transversal intersection of these invariant manifolds. These transversal intersections ensure the existence of a horseshoe, and therefore of chaotic dynamics \cite{Smale65}.

Even if there is a standard theory to analyze the intersection of stable and unstable invariant manifolds, the so called Melnikov Theory (see \cite{Melnikov63} and \cite{GuckenheimerH83} for a more modern exposition), as we will see,  it does not apply to this kind of systems. The reason is that we are in a singular perturbative setting (see Section \ref{sec:RescaledSystem}). Indeed, as A. Neisthadt proved in \cite{Neishtadt84},  in this setting,  when the system is analytic, the invariant manifolds are exponentially close to each other with respect to the parameter $\de$. This implies that there is no hope that classical perturbative techniques, such as Melnikov Theory, will allow to study the possible intersection of these manifolds.

The problem of exponentially small splitting of separatrices was first pointed out by Poincar\'e \cite[\S
19]{Poincare90}  and has drawn considerable attention in the past decades (see  \cite{HolmesMS88, HolmesMS91,DelshamsS92, Fontich93, ChierchiaG94,
Gelfreich94, Fontich95, Sauzin95,DelshamsGJS97, DelshamsS97, Gelfreich97a, Gelfreich97,
Treshev97, GalGM99, Gelfreich00, Sauzin01,OliveSS03,DelshamsGS04, BaldomaF04, BaldomaF05, Baldoma06, Olive06,
GuardiaOS10, BaldomaFGS11, GuardiaS11}). Nevertheless, most of these results cannot be applied to the resonances of nearly integrable systems since they assume an artificial smallness condition for certain terms in the perturbation $h_1$. The only ones that can be applied to this kind of systems are \cite{Treshev97, Gelfreich00, OliveSS03, Olive06, GuardiaOS10}, but they only deal with particular examples.
The present paper is strongly related to \cite{BaldomaFGS11}. In that paper, the authors  study the exponentially small splitting problem for slightly different systems (see the paper for the exact hypotheses). Nevertheless, one can easily see that Hamiltonians \eqref{def:HamResonant} with $h_0(I)=I^2/2$ and $h_1$ independent of $I$ studied in the $0:1$ resonance fit the framework of  \cite{BaldomaFGS11}  and therefore, the results obtained in that paper can be applied to these systems. In the present paper, using the tools developed in  \cite{BaldomaFGS11}, we study the exponentially small splitting problem for more general Hamiltonians  of the form \eqref{def:HamResonant} close to any resonance (in particular allowing general $h_0$ and general dependence on $I$ of $h_1$, see Section \ref{subsec:HypAndResult} for the exact hypotheses). Nevertheless, as happened in  \cite{BaldomaFGS11}, we have to assume that $h_1$ is a trigonometric polynomial in $x$. This hypothesis has been assumed in all the previous works since as far as the author knows, there are not known tools to deal with more general (analytic in $x$) perturbations.

The structure of the paper goes as follows. First, in Section \ref{sec:Results:0:1} we study the splitting of separatrices for Hamiltonians of the form \eqref{def:HamResonant} in the $0:1$ resonance. In this section we state the hypotheses we need to impose to Hamiltonian \eqref{def:HamResonant} (Section \ref{subsec:HypAndResult}) and we state the main results (Section \ref{sec:MainResults}). In Section \ref{sec:Example} we apply the obtained result to a particular (and paradigmatic) example. From the results obtained for the $0:1$ resonance,  in Section \ref{sec:GeneralResonance} we deduce analogous results for all the other $n/m:1$-resonances. Finally in Section \ref{sec:SketchProof} we prove the results stated in the previous sections.

\section{Results for the $0:1$ resonance}\label{sec:Results:0:1}
\subsection{The rescaled Hamiltonian System in the resonance}\label{sec:RescaledSystem}
To facilitate the exposition, in this section we focus on the $0:1$-resonance. Namely we consider a region of the phase space $I\sim I^\ast$ such that $I^\ast$ satisfies that the corresponding frequency vector is  $\omega(I^\ast)=(\pa_I h_0(I^*),1)=(0,1)$. Any other resonance can be reduced to this one after a change of variables as will be seen in Section \ref{sec:GeneralResonance}. After a translation in $I$, one can assume that $\omega(I^\ast)=(0,1)$ takes place at $I^\ast=0$, which implies that $\pa_I h_0(0)=0$. Then $h_0$ is of the form
\begin{equation}\label{def:Original:h_0}
h_0(I)=h_{02}\frac{I^2}{2}+g(I)\,\, \text{ with }\,\,
g(I)=\OO\left(I^3\right).
\end{equation}
We assume that $h_{02}\neq 0$, which guarantees that the resonance is non-degenerate (see Section \ref{subsec:HypAndResult} for the exact hypotheses we need). Then, without loss of generality, we can also assume 
\begin{equation}\label{hyp:PositiveDerivative}
h_{02}>0. 
\end{equation}

The dynamics of the perturbed system around this resonance, is better understood if one performs the rescaling
\begin{equation}\label{def:Rescalings}
I=\sqrt{\frac{\de}{h_{02}}} y\quad\text{ and }\quad \tau=\frac{t}{\sqrt{h_{02}\de}},
\end{equation}
since then we magnify the size of the resonant zone, so that its width becomes of order one with respect to $\de$. To simplify the notation, we take  
\begin{equation}\label{def:NewParameter}
\eps=\sqrt{h_{02}\de}
\end{equation}
 as a new parameter. Then, one  obtains a Hamiltonian System with Hamiltonian
\begin{equation}\label{def:ToyModel:Previ}
H\left(x,y,\frac{t}{\eps};\eps\right)=\frac{y^2}{2}+\frac{1}{\eps^2}G(\eps
y)+V(x)+F\left(x,\frac{t}{\eps}\right)+R\left(x,\eps
y,\frac{t}{\eps};\eps\right),
\end{equation}
where
\begin{align}
\dps G(I)&= h_{02}g\left(\frac{I}{h_{02}}\right) \label{def:Perturb:G}\\
\dps V(x)&=\langle h_1(x,0,\tau;0)\rangle=\frac{1}{2\pi}\int_0^{2\pi}
h_1(x,0,\tau;0)\,d\tau\label{def:PotentialV}\\
\dps F(x,\tau)&= h_1(x,0,\tau;0)-\langle h_1(x,0,\tau;0)\rangle\label{def:Perturb:F}\\
\dps R(x,I,\tau;\eps)&=h_1\left(x,\frac{I}{h_{02}},\tau;\frac{\eps^2}{h_{02}}\right)-h_1(x,0,\tau;0).\label{def:Perturb:R}
\end{align}
Note that  there are terms of different size. The terms $y^2/2$ and $V(x)$ are independent of the new parameter $\eps$ and have become the main terms of the rescaled Hamiltonian. The term $F(x,t/\eps)$ is also of order one with respect to $\eps$ but it is rapidly oscillating and has zero average. Thus, as it is well known from Averaging Theory (see for instance \cite{ArnoldKN88}), this term is in fact causing effects of order $\OO(\eps)$ to the system. More precisely, one can perform an $\eps$-close to the identity symplectic change of variables so that this term becomes of order $\OO(\eps)$. Finally, the terms $\eps^{-2}G(\eps y)$ and $R\left(x,\eps
y,t/\eps;\eps\right)$ are of order $\OO(\eps)$.

Taking into account these different sizes, one can split  Hamiltonian \eqref{def:ToyModel:Previ} as a new unperturbed Hamiltonian plus a new perturbation. Namely, one can rewrite it as follows.
\begin{equation}\label{def:ToyModel}
H\left(x,y,\frac{t}{\eps}\right)=H_0(x,y)+\mu
H_1\left(x,y,\frac{t}{\eps};\eps\right)
\end{equation}
with
\begin{align}
H_0(x,y)&=\frac{y^2}{2}+V(x)\label{def:ToyModel:H0}\\
H_1(x,y,\tau;\eps)&=F(x,\tau)+\frac{1}{\eps^2}G(\eps y)+R(x,\eps y,
\tau;\eps).\label{def:ToyModel:H1}
\end{align}
Note that here $\mu$ is a \emph{fake} parameter since we are interested in $\mu=1$. Nevertheless, we write it to clearly stress that now we have a new unperturbed system $H_0$ (that is, $H$ with $\mu=0$) and a perturbation $\mu H_1$. Moreover, stating the results including the parameter $\mu$ makes them simpler to compare them with the previous ones which usually also include this parameter. Indeed, the first results in the area assumed some smallness condition on this parameter $\mu$ (see \cite{HolmesMS88, DelshamsS92, Gelfreich94,DelshamsS97, Gelfreich97,BaldomaF04, BaldomaF05}). In the case $\mu=1$ there are fewer results. In \cite{Treshev97, Gelfreich00, GuardiaOS10} the authors study particular examples, which in particular satisfy $G=0$ and $R=0$. In \cite{Baldoma06, BaldomaFGS11} the authors study systems of a slightly different form. Nevertheless, one can see that general Hamiltonians of the form \eqref{def:ToyModel} with $G=0$ and $R=0$ fit their framework. Here we generalize these results to Hamiltonians \eqref{def:ToyModel} with general $G$ and $R$. In particular we show that the functions $G$ and $R$ contribute to the main order of the difference between the invariant manifolds and, in particular, if one omits them one obtains a wrong size of the asymptotic order.

The equations associated to the Hamiltonian
\eqref{def:ToyModel} are
\begin{equation}\label{eq:ode:original0}
\left\{\begin{array}{l}\dps \dot x=y+\mu \pa_y
H_1\left(x,y,\frac{t}{\eps};\eps\right)\\
\dps\dot y=-V'(x)-\mu \pa_x
H_1\left(x,y,\frac{t}{\eps};\eps\right).
\end{array}\right.
\end{equation}
From now on, we call  unperturbed system to the system defined by
the Hamiltonian $H_0$ and we refer to $H_1$ as the perturbation.


\subsection{Hypotheses}\label{subsec:HypAndResult}
We devote this section to state the hypotheses needed for the main results concerning the resonance $0:1$. Note that all the hypotheses refer to the integrable Hamiltonian $h_0$ and the first order in $\de$ of the perturbation $h_1$ in \eqref{def:HamResonant}. Namely, the $\OO(\de^2)$ terms in \eqref{def:HamResonant} do not play any role. Some of them are stated using the rescaled hamiltonian \eqref{def:ToyModel}.

\begin{description}
\item[\textbf{HP1}] There exists $r>0$ independent of $\de$, such that the Hamiltonian $h$ in \eqref{def:HamResonant} is analytic for $(x,I,\tau)\in \TT\times B(r)\times\TT$, where $B(r)$ is the ball $B(r)=\{I\in \CC: |I|<r\}$. Moreover, $h_0$ satisfies
\[
 \pa_I^2 h_0(0)> 0.
\]
\item[\textbf{HP2}] The original Hamiltonian is of the form \eqref{def:HamResonant} and $h_1$ is a trigonometric polynomial of degree $M>0$.


\item[\textbf{HP3}] The potential $V$ in \eqref{def:PotentialV} is
 a trigonometric polynomial of the same degree as $h_1$ in \eqref{def:HamResonant}. Namely, $V$ has degree $M$.
\item[\textbf{HP4}] The system associated to $H_0$ has a hyperbolic critical point at $(0,0)$ with
eigenvalues $\{\la,-\la\}$ with $\la>0$. Equivalently, the potential $V$ defined in \eqref{def:PotentialV} satisfies
\[
V(x)=-\frac{\la^2}{2}x^2+\OO\left(x^3\right)\,\quad\text{ as
}x\rightarrow 0.
\]
\item [\textbf{HP5}] The stable and unstable invariant manifolds of the critical point $(0,0)$ of the system associated to  the Hamiltonian $H_0$  coincide along a separatrix.
\end{description}

We denote by $(q_0(u),p_0(u))$ a real-analytic time parametrization of the
separatrix given by Hypothesis \textbf{HP5} with some chosen (fixed) initial condition.
It is well known (see \cite{Fontich95}) that  there exists $a>0$ such that
the parametrization $(q_0(u), p_0(u))$ is analytic in the complex strip $\{ |\Im
u |<a\}$ and not in any wider strip. Next hypothesis deals with the behavior of this parameterization at the boundary of this strip.
\begin{description}
 \item [\textbf{HP6}] There exists a real-analytic time parametrization of the
separatrix $(q_0(u),p_0(u))$ analytic on $\{|\Im u | < a\}$
such that the only singularities of $(q_0(u),p_0(u))$ in the lines $\{\Im u=\pm
a\}$ are at $u=\pm ia$.
\end{description}
Note that the strong assumption in Hypothesis \textbf{HP6} is having only one singularity at each boundary. Once there exists only one, one can just make a time translation to locate it at the imaginary axis. 

Together Hypothesis \textbf{HP2} and \textbf{HP6} imply that  $q_0(u)$ has logarithmic
singularities at $\pm ia$
of the form  $q_0(u)\sim \ln(u\mp ia)$ (where we  take different
branches of the logarithm whether we are close to $+ia$ or $-ia$: we
take $\mathrm{arg}(u-ia)\in (-3\pi/2,\pi/2)$ and
$\mathrm{arg}(u+ia)\in(-\pi/2,3\pi/2)$ respectively). In this case, one can see
that,  if
$u\in\CC$, $|u\mp ia|<\nu$,
\begin{equation}\label{eq:SepartriuAlPolTrig}
\begin{split}
\dps \cos(q_0(u))&=\frac{\wh C^1_\pm}{(u\mp
ia)^{2/M}}\left(1+\OO\left((u\mp ia)^{2/M}\right)\right)\\
\dps \sin(q_0(u))&=\frac{\wh C^2_\pm}{(u\mp
ia)^{2/M}}\left(1+\OO\left((u\mp ia)^{2/M}\right)\right)\\
\dps p_0(u)&=\frac{C_\pm}{(u\mp ia)}\left(1+\OO\left((u\mp
ia)^{2/M}\right)\right)
\end{split}
\end{equation}
with  $\mathrm{arg}(u-ia)\in (-3\pi/2,\pi/2)$ and
$\mathrm{arg}(u+ia)\in(-\pi/2,3\pi/2)$ if we are dealing with the
singularity $+ia$ or $-ia$ respectively (see Section 2.1.3 of \cite{BaldomaFGS11}). We also have that 
\begin{equation*}
C_+ = \overline{C_-}= \pm i\frac{2}{M} 
\end{equation*}

\begin{remark}\label{remark:NonHetero}
Hypothesis \textbf{HP1} can be weakened to  $\pa_I^2 h_0(0)\neq 0$ and both the positive and negative case can be proved analogously. To simplify the exposition we assume  $\pa_I^2 h_0(0)> 0$ (the other case can be deduced from it just reversing time). Note that this condition has already been stated in  \eqref{hyp:PositiveDerivative} and has been used for the rescaling and obtention of the model close to the resonance \eqref{def:HamResonant}.
\end{remark}
\begin{remark} The critical point $(0,0)$ has invariant manifolds thanks to Hypothesis \textbf{HP4}. Therefore, the actual assumption in Hypothesis \textbf{HP5} is that they coincide along a separatrix. It is easy to see that the energy levels of $H_0$ are compact sets and therefore, the only possibilities is that the invariant manifolds of $(0,0)$ coincide with either the invariant manifolds of the same or another critical point. Therefore, \textbf{HP5} rules out the case of heteroclinic connections. The existence of heteroclinic connections is a non-generic phenomena, since after a small perturbation one can make that all the critical points of $H_0$ belong to different energy levels and thus impossible to be connected by heteroclinic connections (see Section \ref{subsec:Genericity}).
\end{remark}

\begin{remark}\label{remark:SeparatrixType}
 The separatrices, whose existence is guaranteed by Hypothesis \textbf{HP5}  can be of two different topological type. Indeed, note that the hyperbolic critical points of the system associated to Hamiltonian \eqref{def:ToyModel:H0}  are of the form $(x^\ast,0)$, where $x^\ast$ are non-degenerate maxima of $V$. Then, if the maximum is global, the associated separatrices are a graph over the variable $x$, as happens for the pendulum, which has only one maximum (see Figure \ref{fig:Pendol}). If the maximum is local but not global, the associated  separatrices form a figure eight as can be seen in Figure \ref{fig:FigureEight}. Therefore, if $V$ has more than one maximum there always coexist graph-like and figure eight separatrices associated to different critical points. The results obtained in this paper can be applied to both types of separatrices as long as Hypotheses \textbf{HP1}-\textbf{HP6} are satisfied.
\end{remark}

\begin{figure}[h]
\begin{center}
\psfrag{xx}{$x$}\psfrag{xxxxx}{$y$}
\includegraphics[height=6cm]{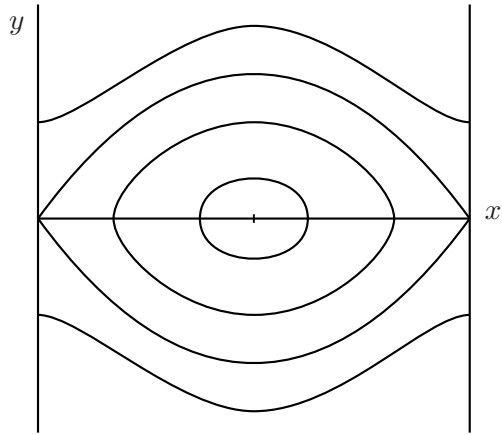}
\end{center}
\caption{\figlabel{fig:Pendol} Phase portrait of the system associated to Hamiltonian \eqref{def:ToyModel:H0} when the potential $V$ has only one local (and thus global) maximum, which after a translation in $x$ can be placed at the origin. Note that this maximum has to be non-degenerate to ensure that the critical point is hyperbolic. See Figure \ref{fig:FigureEight} for the phase portrait when $V$ has more local maxima.}
\end{figure}

\begin{figure}[ht]
\begin{center}
\psfrag{x}{$x$}\psfrag{y}{$y$}\psfrag{0}{$0$}\psfrag{P}{$2\pi$}
\includegraphics[width=10cm]{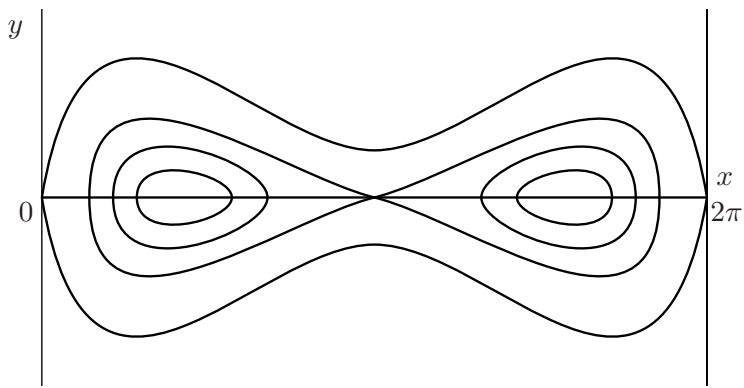}
\end{center}
\caption{\figlabel{fig:FigureEight}  Phase portrait of the system associated to Hamiltonian \eqref{def:ToyModel:H0} when the potential $V$ has one global and one local maximum. The global one corresponds to the critical point whose separatrices are a graph over $x$. The local (and non-global) critical point is the one having figure eight separatrices. If $V$ has more local maxima, one obtains more nested figure eight separatrices in the phase space}
\end{figure}

\subsection{Main results for the $0:1$ resonance}\label{sec:MainResults}
We devote this section to state the main results about the $0:1$ resonance. First we state  Theorems \ref{th:MainPO} and \ref{th:Main}, which give the main results referred to the rescaled Hamiltonian \eqref{def:ToyModel} and then, from them, we deduce Theorem \ref{th:Main:Rescaled:(0,1)}, which states  the main results for the resonance $0:1$ in the original variables. That is, we give it in terms of  Hamiltonian \eqref{def:HamResonant}. In Section \ref{sec:GeneralResonance} we deduce from  Theorems \ref{th:MainPO} and \ref{th:Main} analogous results for general resonances $n/m:1$.

By Hypothesis \textbf{HP4}, the system associated to the Hamiltonian $H_0$ in  \eqref{def:ToyModel:H0}  has a  hyperbolic critical  point at the origin. Next theorem ensures that the  hyperbolic critical point of
the unperturbed system becomes a hyperbolic periodic orbit for the full system \eqref{eq:ode:original0}, which is 
close to the origin.

\begin{theorem}\label{th:MainPO}
Let us  assume Hypotheses \textbf{HP1} and \textbf{HP4}. Then, given any fixed value $\mu>0$,
there exists $\eps_0>0$ such that if
$\eps\in(0,\eps_0)$, the system associated to \eqref{eq:ode:original0} has a hyperbolic
periodic orbit $(x_p(t/\eps),y_p(t/\eps))$ which satisfies that, for
$\tau\in\RR$,
\[
\begin{split}
\left
|x_p\left(\tau\right)\right|&\leq
K |\mu|\eps^{2}\\
\left|y_p\left(\tau
\right)\right|&\leq
K |\mu|\eps,
\end{split}
\]
for a constant $K>0$ independent of $\eps$.
\end{theorem}
The proof of this theorem follows the same lines as the proof of Theorem 2.2 in \cite{BaldomaFGS11} (see also \cite{Fontich93,Fontich95,DelshamsS97}).

The next step is to study the stable and unstable invariant
manifolds of the periodic orbit $(x_p,y_p)$. In the unperturbed case
(that is $\mu=0$) we know that they coincide along the separatrix
$(q_0,p_0)$ given in \textbf{HP6}. We will see that when $\mu\neq 0$ they tipically split and their transversality is exponentially small.

To measure the splitting of the invariant manifolds we consider the
$2\pi\eps$-Poincar\'e map $P_{t_0}$ in a transversal section
$\Sigma_{t_0}=\left\{(x,y,t_0); (x,y)\in\RR^2\right\}$. This
Poincar\'e map has a hyperbolic  fixed point
$(x_p(t_0/\eps),y_p(t_0/\eps))$, whose stable and unstable invariant manifolds are curves.

As $P_{t_0}$ is an area
preserving map, we measure the splitting giving an asymptotic
formula for the area of the lobes generated by these curves between
two transversal homoclinic points. Moreover, by the area preserving
character of $P_{t_0}$, the area $\AAA$ of these lobes does not
depend on the choice of the homoclinic points. Other quantities measuring
the splitting, as the distance along a transversal section to the unperturbed
separatrix, or the angle between these curves at an homoclinic point, can be
easily derived from our work.

\begin{theorem}\label{th:Main}
Let us assume Hypotheses \textbf{HP1}-\textbf{HP6}. Then, given any fixed $\mu$,
there exists $\eps_0>0$ such that if $\eps\in (0 ,\eps_0)$, the area of the
lobes between the invariant manifolds of the periodic orbit given in Theorem
\ref{th:MainPO}  is given by
\begin{equation}\label{def:formulaArea:singular}
\AAA=4|\mu|\eps^{-1}e^{-{\dps\tfrac{a}{\eps}}+\mu\Im
b\ln\frac{1}{\eps}}\left(\left| f
\left(\mu\right)e^{iC\left(\mu\right)}\right|+\OO\left(\frac{1}{|\ln\eps|}
\right)\right),
\end{equation}
where  $f(\mu)$ is an entire analytic function, $C(\mu)$ is
an entire analytic function such
that $C(\mu)=\OO(\mu)$ and $b\in \CC$ is a constant defined as 
\begin{equation}\label{def:Constantb}
b=2g_3C_+,
\end{equation}
where $C_+$  has been defined in \eqref{eq:SepartriuAlPolTrig} and $g_3$ is the degree 3 Taylor coefficient of the function $G$ in \eqref{def:Perturb:G}, namely, $G(I)=g_3 I^3+\OO(I^4)$.
\end{theorem}

The proof of this theorem is given in Section \ref{sec:SketchProof}. Recall that we were primarily interested in the case $\mu=1$, in this case the asymptotic formula for the area becomes 
\begin{equation}\label{def:FormulaAreaMu1}
 \AAA=4\eps^{-1-\Im b}e^{-{\dps\tfrac{a}{\eps}}}\left(\left| f
\left(1\right)e^{iC\left(1\right)}\right|+\OO\left(\frac{1}{|\ln\eps|}
\right)\right).
\end{equation}

Theorem \ref{th:Main} can be rephrased to give an analogous result for the system associated to the original Hamiltonian $h$ in \eqref{def:HamResonant}.

\begin{theorem}\label{th:Main:Rescaled:(0,1)}
Let us assume Hypotheses \textbf{HP1}-\textbf{HP6}. Then, 
there exists $\de_0>0$ such that if $\de\in (0 ,\de_0)$, the system associated to Hamiltonian \eqref{def:HamResonant} has a hyperbolic periodic orbit $(x_p(\tau), I_p(\tau))$ which satisfies
\[
 \begin{split}
  \left|x_p(\tau)\right|&\leq K \de\\
\left|I_p(\tau)\right|&\leq K \de.
 \end{split}
\]
Moreover, the area of the lobes between its invariant manifolds is given by
\begin{equation}\label{def:FormulaLobeOriginal}
 \wt \AAA=4({h_{02}})^{-1-{\dps\tfrac{\Im b}{2}}}\de^{-{\dps\tfrac{\Im b}{2}}}e^{-{\dps\tfrac{a}{\sqrt{h_{02}\de}}}}\left(\left| \Theta\right|+\OO\left(\frac{1}{|\ln\de|}
\right)\right),
\end{equation}
where $h_{02}=\pa_I^2 h_0(0)$ (which is positive by Hypothesis \textbf{HP1}),  $\Theta = f(1) e^{iC(1)}$ and $b$, $f(1)$ and $C(1)$ are the constants given by Theorem \ref{th:Main}.
\end{theorem}

Note that to deduce this theorem from Theorem \ref{th:Main} it is enough to undo the changes \eqref{def:Rescalings} and \eqref{def:NewParameter}.

\subsubsection{Some remarks}

In the first order of formulas \eqref{def:FormulaAreaMu1} and \eqref{def:FormulaLobeOriginal}   there appear three constants $f(1)$, $C(1)$ and $b$. We devote this section to make some comments about them.
\begin{itemize}
\item  The constants $f(1)$ and  $C(1)$ play a similar role. Nevertheless, the functions $f(\mu)$ and $C(\mu)$ have a significantly different origin, which can be seen in Section \ref{sec:SketchProof}. Let us just say here that $C(\mu)$ can be obtained by classical perturbative methods, in the sense that it  can be computed in terms of the first orders in $\eps$ of the funtions $R(x,\eps y,\tau;\eps)$ and $\eps^{-2}G(\eps y)$. Instead, $f(\mu)$ has a completely different origin. It comes from the study of the \emph{inner equation} associated to system \eqref{eq:ode:original0}, which is an equation independent of $\eps$ and was studied in \cite{Baldoma06} (see Section \ref{sec:SketchProof} for its definition). In Section \ref{sec:SketchProof}, we will see that it depends on the full jet in $I$ of the funtions $G(I)$ and $R(x,I,\tau;\eps)$. Namely, it depends on all the orders in $\eps$ of  $R(x,\eps y,\tau;\eps)$ and $\eps^{-2}G(\eps y)$. Note that this implies that any finite order truncation in $\eps$ of the Hamiltonian $H$ does not predict correctly the area of the lobes.
\item To ensure that the area of the lobe is indeed given asymptotically by  formula  \eqref{def:FormulaAreaMu1}, one has to check that $f(1)\neq0$. With the techniques used in this work, as far as the authors know, there is not any analytical way to check this condition. Nevertheless, one can use the mentioned inner equation to verify this condition numerically \cite{SimoV10}.
\item The constant $b$  affects the polynomial term in the formulas for the area of the lobes  \eqref{def:FormulaAreaMu1} and \eqref{def:FormulaLobeOriginal}. Generically satisfies $b\neq 0$, since generically $g_3\neq 0$. Nevertheless, it vanishes when $g_3=0$. In previous works (see \cite{Treshev97, Olive06}), the models considered satisfied $G=0$ and thus this coefficient did not appear in the formulas provided by  those works. This term was first detected in \cite{Baldoma06, BaldomaFGS11}. Note that, as it happens for the function $C(\mu)$, it can be obtained by the first order in $G$, that is, one does not need the full jet of this function to compute it, as happens for $f(\mu)$.
\end{itemize}

\subsubsection{Comparison with Melnikov Theory}\label{sec:Melnikov}
If one considers the system associated to Hamiltonian \eqref{def:ToyModel}, one can apply the classical Melnikov Theory, which is just to make a perturbative approach taking  $\mu$ as  a small parameter. Recall that we are interested in the case $\mu=1$ and not  in arbitrarily small $\mu$. Nevertheless, this was the first approach to study the phenomenon of exponentially small splitting of separatrices since turns out to be the simplest case. Indeed, using this method and taking $\mu$ small enough, one can see that the area of the lobes is given asymptotically in $\mu$ as
\begin{equation}\label{def:MelnikovPrediction}
 \AAA=4\eps\ii e^{-{\dps\tfrac{a}{\eps}}}\left(\left| f_0
\right|\mu+\OO\left(\mu^2\right)\right),
\end{equation}
for certain constant $f_0\in\CC$.

Thus, if one wants to assure the validity of the first order given by this formula, it has to be assumed that $\mu$ is exponentially small with respect to $\eps$. Nevertheless, it is a well known fact that Melnikov gives the correct first order for a wider range in $\mu$ (see \cite{DelshamsS97, BaldomaFGS11}). Using Theorem \ref{th:Main} one can see for which range of $\mu$ Melnikov indeed gives the correct first order and for which range it does not. In \cite{Baldoma06}, it is seen that the function $f(\mu)$ appearing in Theorem \ref{th:Main} and the constant $f_0$ involved in the Melnikov prediction \eqref{def:MelnikovPrediction} are related by
\[
 f(\mu)=f_0+\OO\left(\mu\right).
\]
Thus, comparing formulas \eqref{def:formulaArea:singular} and \eqref{def:MelnikovPrediction}, one can see that in the general case $b\neq 0$, the Melnikov prediction is correct provided
\[
 \mu\ll \frac{1}{|\ln\eps|}.
\]
In the non-generic case $b=0$, which takes place for instance when $G=0$,  it is enough to require $\mu\ll1$.

In the resonances of nearly integrable Hamiltonian Systems, the parameter $\mu$ satisfies $\mu=1$ and therefore Melnikov does not predict correctly the area of the lobes. Indeed, comparing formulas \eqref{def:formulaArea:singular} and \eqref{def:MelnikovPrediction}, one can conclude that  Melnikov only predicts correctly the coefficient in the exponential, but fails to predict correctly both the polynomial term and the constant coefficient.

\subsubsection{On the genericity of the hypotheses and the results}\label{subsec:Genericity}
We devote this section to make some remarks on the genericity of Hypotheses \textbf{HP1-HP6} stated in Section \ref{subsec:HypAndResult}, assumed in Theorem \ref{th:Main} and \ref{th:Main:Rescaled:(0,1)}, and about the genericity of the results stated in these theorems.

Let us consider a constant $r>0$. Then, it is clear that \textbf{HP1}, which is the only condition needed for the integrable system $h_0$ in \eqref{def:HamResonant}, is satisfied for a generic $h_0$ in 
\begin{equation}\label{def:Spaceh0}
 \XX=\left\{f:\{|I|\leq r\}\subset\CC\rightarrow \CC: \text{real-analytic}, \|f\|_{\infty}<\infty\right\},
\end{equation}
where $\|\cdot\|_\infty$ denotes the classical supremmum norm.

Concerning the hypotheses imposed on the perturbation $h_1$ in \eqref{def:HamResonant}, let us fix also constants $\eps_0>0$ and $\sigma>0$. Then, one would like to obtain their genericity in the space
\[
 \YY=\left\{f:\TT_\sigma\times\{|I|\leq r\}\times\TT_\sigma\times \{|\eps|\leq \eps_0\}\rightarrow \CC: \text{real-analytic}, \|f\|_{\infty}<\infty\right\},
\]
where $\TT_\sigma=\{x\in \CC/\ZZ: |\Im x|\leq \sigma\}$. Nevertheless, functions $h_1\in\YY$ satisfying  Hypothesis \textbf{HP2} are clearly non generic since they are trigonometric polynomials in $x$. However, one can ask about the genericity in the smaller space of perturbations which are a trigonometric polynomial in $x$ of degree $M>0$ or less. Thus, we consider the space
\begin{equation}\label{def:Spaceh1}
\begin{split}
 \YY_M=\bigg\{f=\sum_{k=-M}^M f^{[k]}(I,\tau;\eps)e^{ikx}:\TT_\sigma\times\{|I|\leq r\}\times&\TT_\sigma\times \{|\eps|\leq \eps_0\}\rightarrow \CC;\\
 &\text{real-analytic}, \|f\|_{\infty}<\infty\bigg\},
\end{split}
\end{equation}
which is a Banach space  contained in $\YY$. By definition, \textbf{HP2} is generic in this space. Taking into account the definition of $V$ in terms of $h_1$ in \eqref{def:PotentialV}, one can easily see that \textbf{HP3} is also generic in $\YY_M$. The same happens for \textbf{HP4} since, as explained in Remark \ref{remark:SeparatrixType}, the hyperbolic critical points of  the system $H_0$ in \eqref{def:ToyModel:H0}  correspond to non-degenerate maxima of the potential $V$ and it is a well known fact that for a generic $V$ all its extrema are non-degenerate. Moreover, in Remark \ref{remark:NonHetero}, we have already explained that having a homoclinic connection for a hyperbolic fixed point is also generic.  Thus we can conclude that Hypotheses \textbf{HP2}-\textbf{HP5} are generic.

It only remains to check \textbf{HP6}. It is clear that this hypothesis is open in $\YY_M$. Let us check its genericity. As we have seen in Remark \ref{remark:SeparatrixType} there are two types of separatrices: the ones that are graphs with respect to $x$ and the ones that form a figure eight (see Figures \ref{fig:Pendol} and \ref{fig:FigureEight}). Let us fix a concrete separatrix and check the genericity of Hypothesis \textbf{HP6}.

In the figure eight case, the separatrices are reversible with respect to the involution $\Phi(x,y)=(x,-y)$. Then, Hypothesis \textbf{HP6} is open but not generic. Indeed, if for a concrete potential, a figure eight separatrix has singularities at points $u=\pm \alpha\pm  ai$ with $\alpha\neq 0$, this structure cannot be broken by any small perturbation of the potential due to the reversibility.

Namely,  for this kind of separatrices in $\YY_M$ there  is an open set where Hypothesis \textbf{HP6} is satisfied and a different open set where the separatrix has two singularities at each connected component of its strip of definition.

For the graph  separatrices, as a far as the author knows, there are no results in the literature about the genericity of Hypothesis \textbf{HP6}. Nevertheless, since in this case the separatrices generically do not have any symmetry, one would expect that  if a separatrix  has singularities at $u=\pm \alpha\pm  ai$, one can modify slightly the potential $V$ so that the singularities have different imaginary part. This fact, prompt us to state the following conjecture.

\begin{conjecture}\label{conj:singularities}
Let us consider a Hamiltonian of the form 
\[
 H_0(x,y)=\frac{y^2}{2}+V(x).
\]
Then, for a generic trigonometric polynomial $V$, the associated Hamiltonian system has a hyperbolic critical point, such that its invariant manifolds coincide forming a separatrix which is a graph with respect to $x$ and satisfies Hypothesis \textbf{HP6}.
\end{conjecture}

Recall that for a generic potential there always exists a hyperbolic critical point having graph separatrices, since they correspond to a non-degenerate global maximum of $V$. Thus, if this conjecture is true, the assumed hypotheses are satisfied by a generic $h_0\in\XX$ (see  \eqref{def:Spaceh0}) and a generic $h_1\in\YY_M$  (see \eqref{def:Spaceh1}). Nevertheless,  having generic hypotheses on the Hamiltonians $h_0$ and $h_1$ does not imply directly that we can deduce from  our results that generically the invariant manifolds split transversally and that the area of the associated lobes is given by the asymptotic formula \eqref{def:FormulaLobeOriginal} with non-vanishing first order. Indeed, to obtain this result one would need to prove that $f(1)\neq 0$ for a generic set of Hamiltonians. 

Let us show that assuming Conjecture \ref{conj:singularities} and adding an extra generic hypothesis this is the case (under a suitable definition of genericity). Indeed, as stated in Theorem \ref{th:Main}, the function $f(\mu)$ is analytic. Moreover, recall that $f(\mu)=f_0+\OO(\mu)$ where  $f_0$ is the constant given by Melnikov Theory (see Section \ref{sec:Melnikov}). This constant $f_0$ comes from certain residue that generically does not vanish, which implies that generically $f_0\neq 0$ (see \cite{DelshamsS97, BaldomaF04} for a more detailed explanation). Then, since $f$ is analytic, under this extra hypothesis $f$ only vanishes for a discrete set  of values of $\mu$ and thus for a generic set of Hamiltonians and a generic value of $\mu$, we have that $f(\mu)\neq 0$. 
Now we show that this fact is sufficient to show that generically $f(1)\neq 0$. The openness of this condition is clear. Therefore, it is enough to check that it is dense. Let us consider $h_0$ and $h_1$ satisfying all the stated hypotheses (and that $f_0$ satisfies $f_0\neq 0$) such that  $f(1)$ vanishes. Since $f(\mu)$ only vanishes for a finite set of points, there exists $\mu^*$  arbitrarily close to 1 such that $f(\mu^*)\neq0$. Notice that modifying slightly $\mu$ means modifying slightly $H_1$, which is equivalent to modifying slightly $h_0$ and $h_1$ in \eqref{def:HamResonant}. Thus, for generic pair $(h_0,h_1)\in\XX\times\YY_M$ (with respect to the product topology), we have that $f(1)\neq 0$.

Therefore, we can conclude that, if we fix $M>0$ and  assume Conjecture \ref{conj:singularities}, then for a generic  pair $(h_0,h_1)\in\XX\times\YY_M$ of Hamiltonians, there exists $\eps_0>0$ such that for any $\eps\in (0,\eps_0)$ the Hamiltonian \eqref{def:HamResonant} has a hyperbolic periodic orbit at the resonance whose invariant maniolds are graphs,  intersect transversally and the area of the lobes is given by the asymptotic formula \eqref{def:FormulaLobeOriginal} of Theorem 
\ref{th:Main:Rescaled:(0,1)}, whose first order does not vanish. 

Finally, let us point out that, if one takes into account the rescaled formula \eqref{def:formulaArea:singular}, we can ensure that for a generic set of Hamiltonians (in the sense just explained),  Melnikov Theory does the wrong prediction for the area of the lobes both in the polynomial and constant terms (see Section \ref{sec:Melnikov}).

\subsubsection{The parabolic case}
 The generic Hypothesis \textbf{HP4} assumes the that the Hamiltonian $H_0$ in \eqref{def:ToyModel:H0} has a hyperbolic critical point at the origin. Nevertheless, the results obtained in Theorem \ref{th:Main} can be generalized,  under suitable hypotheses, to Hamiltonian Systems which  have a parabolic point. Indeed,  let us assume that $H_0$ has a parabolic point at the origin which has invariant manifolds that coincide along a separatrix, as happens, for instance, for
\[
 H_0(x,y)=\frac{y^2}{2}-\left(\cos x-1\right)^2.
\]
 Having a parabolic critical point at the origin, is equivalent to ask that the potential  $V(x)$ (see \eqref{def:PotentialV}) satisfies
\[
 V(x)=v_k x^{k}+\OO\left(x^{k+1}\right)\,\,\,\text{ as }x\rightarrow 0
\]
for certain $k\geq 3$ and $v_k\in\RR$. Then,  to avoid bifurcations, the extra hypotheses that one has to assume are to guarantee that, when we add the perturbation, the critical point preserves its parabolic character and that its invariant manifolds are preserved. To this end one has to impose the  following two conditions. The first one is that the function $F$ in \eqref{def:Perturb:F} satisfies
\[
 F(x,\tau)=\alpha(\tau)x^m+\OO\left(x^{m+1}\right)\,\,\,\text{ as }x\rightarrow 0
\]
with $2m-2\geq k$.  The second condition is the analogous one for $R$ in \eqref{def:Perturb:R}. Indeed, if we Taylor-expand $R$ with respect to $I$ around $I=0$ we get
\[
 R(x,I,\tau;\eps)=R_1(x,\tau;\eps)I+\OO\left(I^2\right)
\]
Then, the second additional hypothesis is that $R_1$  satisfies
\[
R_1(x,\tau;\eps)=\beta(\tau;\eps)x^n+\OO\left(x^{n+1}\right)\,\,\,\text{ as }x\rightarrow 0
\]
 with $n\geq k/2$. Under these two extra hypotheses, Theorem \ref{th:Main} is also valid for parabolic points. The proof of this fact was done for slightly different Hamiltonian Systems in \cite{BaldomaFGS11} (see also \cite{BaldomaF04}), and it can be easily adapted to our setting. Note that these two extra hypotheses are satisfied in the natural case that $F$ and $R$ are of the same order as $V$, that is if $F=\OO( x^k)$ and $R=\OO( x^k)$. However, the hypotheses we need to require are considerably weaker.

\subsection{An example}\label{sec:Example}
In this section we apply Theorems \ref{th:Main} and \ref{th:Main:Rescaled:(0,1)} to a particular example. We consider the Hamiltonian 
\begin{equation}\label{def:Example:Ham}
h(x,I,\tau;\de)=\frac{I^2}{2}+\eta I^3+\de \left(1+\alpha I\right)(1+\sin\tau)(\cos x-1),\,\,\,\eta,\alpha\in\RR,
\end{equation}
and we study its splitting of separatrices for the resonance $0:1$, which takes place at $I=0$. We have chosen this model for two reasons. Firstly, (most of) the constants involved in the splitting formula can be computed explicitly and thus we can analyze it in great detail. Secondly, it  shows all the modifications that the splitting formula of Theorem \ref{th:Main} presents in comparison to the Melnikov prediction. Moreover, it can be compared with the results in \cite{Treshev97}, which considers the same model with $\eta=0$ and $\alpha=0$.

Note that for this model $h_{02}=1$, and thus the rescaling is simply $y=\sqrt{\de} I$ and $\eps=\sqrt{\de}$. Then, it is easy to see that the rescaled model is given by
\begin{equation}\label{def:Example:Ham:rescaled}
 H\left(x,y,\frac{t}{\eps}\right)=  H_0(x,y)+H_1\left(x,y,\frac{t}{\eps}\right)
\end{equation}
with
\[
 \begin{split}
  H_0(x,y)&=\frac{y^2}{2}+\cos x-1\\
  H_1(x,y,\tau;\eps)&=\sin \tau (\cos x-1)+\alpha \eps y (1+\sin\tau)(\cos x-1)+\eta\eps y^3
 \end{split}
\]
 Thus, the  functions $V$, $G$ and $F$ and $R$ associated to this system (see \eqref{def:PotentialV}, \eqref{def:Perturb:G}, \eqref{def:Perturb:F} and \eqref{def:Perturb:R}) are given by
\[
 \begin{split}
  V(x)=&\cos x-1\\
  G(I)=& \eta I^3\\
  F(x,\tau)&=\sin \tau(\cos x-1)\\
  R(I,x,\tau)&=\alpha I(1+\sin\tau)(\cos x-1)
 \end{split}
\]
The Hamiltonian $H_0$ is the Hamiltonian of the classical pendulum and therefore   Hamiltonian \eqref{def:Example:Ham} satisfies Hypothesis \textbf{HP4} and \textbf{HP5}. The upper separatrix of the pendulum can be parameterized as
\[
 \ga(u)=\left(4\arctan\left(e^u\right),\frac{2}{\cosh u}\right)
\]
and its singularities are located at $u=i\pi/2+k\pi i$, $k\in\ZZ$. Therefore Hypothesis \textbf{HP6} is also satisfied. It is easy to see that Hamiltonian \eqref{def:Example:Ham} also satisfies Hypotheses \textbf{HP1}-\textbf{HP3} and thus one can apply Theorem \ref{th:Main}. Note also that the periodic orbit at $(x,y)=(0,0)$ is preserved when one adds the perturbation and, therefore, Theorem \ref{th:MainPO} is not needed. 

We have that $g_3=\eta$ and $C_+=-2i$ and therefore, using formula \eqref{def:Constantb}, the constant $b$ is given by $b=-4\eta i$. Then, the asymptotic formula for the splitting for system \eqref{def:Example:Ham:rescaled} reads 
\begin{equation}\label{def:Example:Splitting}
 \AAA=4\eps^{-1+4\eta}e^{-{\dps\tfrac{\pi}{2\eps}}}\left(\left| f
\left(1\right)e^{iC\left(1\right)}\right|+\OO\left(\frac{1}{|\ln\eps|}
\right)\right). 
\end{equation}
As we have already explained, the constant $f(1)$ cannot be computed with our methods and therefore, we cannot ensure whether it vanishes or not. Even if it is not written explicitly, both $f(1)$ and $C(1)$ depend analytically on the parameters $\eta$ and $\alpha$.

If we take $\eta=0$ and $\alpha=0$ we recover from \eqref{def:Example:Splitting} the formula given in  \cite{Treshev97}. As we have said only when $\eta=0$, the constant $b$ vanishes and then, the degree of the polynomial term coincides with the Melnikov prediction.  Note also that one could add higher orders in $I$ both in $h_0$ and $h_1$. Then, in the splitting formula \eqref{def:Example:Splitting} only would change $f(1)$ and $C(1)$ but $b$ would remain the same (see \eqref{def:Constantb}).

Undoing the rescalings, one can deduce from \eqref{def:Example:Splitting} the formula for the original system \eqref{def:Example:Ham}, which reads
\[
 \wt \AAA=4\de^{2\eta}e^{-{\dps\tfrac{\pi}{2\sqrt{\de}}}}\left(\left| f
\left(1\right)e^{iC\left(1\right)}\right|+\OO\left(\frac{1}{|\ln\de|}
\right)\right).
\]

\section{Results for a general resonance}\label{sec:GeneralResonance}

The results stated in Theorem \ref{th:Main:Rescaled:(0,1)} are focused on the resonance with frequency vector $\omega=(0,1)$. Nevertheless, as we will show in this section, from Theorem \ref{th:Main}, one can deduce results for any other resonance. Indeed, we will reduce any resonance to the resonance $0:1$.

As it is well known, a resonance takes place at $I=I^\ast$ if $\pa_I h_0(I^\ast)\in \QQ$. Let us assume, without loss of generality that the resonance takes place at $I^\ast=0$ and then $\pa_I h_0(0)=n/m$ where $n$ and $m$ are coprime integers. We consider a Hamiltonian $h_0+\de h_1$ of the form \eqref{def:HamResonant} and we show that there exists a change of variables such that the new Hamiltonian $\wt h_0+\de \wt h_1$ satisfies $\pa_I \wt h_0(0)=0$. Namely, we  transform the resonance $m/n:1$ into $0:1$. Indeed, this change is simply given by
\begin{equation}\label{def:ChangeAnyResonance}
  \wt x=x-\frac{n}{m}\tau,
\end{equation}
and the transformed Hamiltonian reads
\[
 \wt h(\wt x, I, \wt \tau;\de)=\wt h_0( I)+\de \wt h_1\left(\wt x, I, \wt \tau;\de\right)
\]
   with
\[
\begin{split}
 \wt h_0( I)&=h_0(I)-\frac{n}{m}I\\
\wt h_1(\wt x, I, \tau;\de)&=h_1\left(\wt x-\frac{n}{m}\tau, I, \tau;\de\right).
\end{split}
\]
Note that with this change the perturbation becomes $2\pi m$-periodic in $\tau$. Then, one can apply the rescaling in \eqref{def:Rescalings}. Nevertheless, to proceed as for the $0:1$-resonance, we want that after the rescaling the perturbation $H_1(\wt x,y,t/\eps;\eps)$ is  $2\pi\eps$ periodic in time (in the variable $t$). To this end, we just need to define as a new parameter
\[
 \eps=m\sqrt{h_{02}\de}
\]
instead of \eqref{def:NewParameter}.

Note that $\wt h_{02}=h_{02}$ since the change \eqref{def:ChangeAnyResonance} only modifies the first coefficient of the Taylor series of $h_0$. 

After the rescaling we get a system of the form \eqref{def:ToyModel} with the modified functions $V$, $F$, $G$ and $R$ defined as
\[
 \begin{split}
  V\left(\wt x\right)&=\frac{1}{2\pi}\int_0^{2\pi} h_1 \left(\wt x +\frac{n}{m}\tau, 0,\tau\right)d\tau\\
F\left(\wt x,\tau\right)&=h_1 \left(\wt x +\frac{n}{m}\tau, 0,\tau\right)-V\left(\wt x\right)\\
R\left(\wt x,I,\tau;\eps\right)&=h_1\left(\wt x+\frac{n}{m}\tau, \frac{I}{h_{02} m},\tau;\frac{\eps^2}{h_{02}m^2}\right)-h_1 \left(\wt x+\frac{n}{m}\tau, 0,\tau;0\right)\\
G(I)&=m^2h_{02}g\left(\frac{I}{m h_{02}}\right).
 \end{split}
\]

 Then, from Theorem \ref{th:Main} we can deduce an asymptotic formula for the splitting of separatrices for any resonance. Nevertheless, note that now the parameter $\eps$ is given by $\eps=m\sqrt{h_{02}\de}$ and therefore the necessary smallness condition on the perturbation, $\eps\in (0,\eps_0)$, depends on $m$, namely on the resonance: the higher the resonance is, the smaller $\de$ is needed. 

\begin{theorem}\label{th:AnyResonance}
Let us assume Hypotheses \textbf{HP1}-\textbf{HP6} and let us fix a rational number $n/m$, where $n$ and $m$ are coprime integers. Then, 
there exists $\de_0>0$ such that if $\de\in (0 ,\de_0)$, the system associated to Hamiltonian \eqref{def:HamResonant} has a hyperbolic periodic orbit $(x_p(\tau), I_p(\tau))$ which satisfies
\[
 \begin{split}
  \left|x_p(\tau)-\frac{n}{m}\tau\right|&\leq K m^2\de\\
\left|I_p(\tau)\right|&\leq K \de.
 \end{split}
\]
Moreover, the area of the lobes between the invariant manifolds of this hyperbolic orbit is given by
\[
 \AAA=4h_{02}^{-1-{\dps\tfrac{\Im b}{2}}}\de^{-{\dps\tfrac{\Im b}{2}}}e^{-{\dps\tfrac{a}{m\sqrt{h_{02}\de}}}}\left(\left| f
\left(1\right)e^{iC\left(1\right)}\right|+\OO\left(\frac{1}{|\ln(m\de^{1/2})|}
\right)\right),
\]
where $b, f(1), C(1)$ are the constants given by Theorem \ref{th:Main} and $h_{02}=\pa_I^2 h_0(0)$, which is positive by Hypothesis \textbf{HP1}.
\end{theorem}

\section{Description of the proof of Theorem \ref{th:Main}}\label{sec:SketchProof}
The proof of Theorem \ref{th:Main} follows the same lines as the proof of Theorem 2.7 in \cite{BaldomaFGS11}. In this section, we describe the proof focusing on the main differences from \cite{BaldomaFGS11}. 

\subsection{Basic notations}\label{sec:basicnotation}
First, we introduce some basic notations which is used in this section. 
We denote by $\TT=\RR/(2\pi\ZZ)$ the real 1-dimensional torus and by
\[
\TT_\sigma=\left\{\tau\in \CC/(2\pi\ZZ); |\Im \tau|<\sigma\right\},
\]
with $\sigma>0$, the torus with a complex strip.

Given a function $h:D\times\TT_\sigma\rightarrow \CC$,
where $D\subset\CC$ is an open set, we denote its Fourier
series by
\[
h(u,\tau)=\sum_{k\in\ZZ}h^{[k]}(u)e^{ik\tau}
\]
and its average by
\[
\langle h\rangle (u)
=h^{[0]}(u)=\frac{1}{2\pi}\int_0^{2\pi}h(u,\tau)\,d\tau .
\]


We  fix also some notation for the functions $F$, $G$ and $R$ defined in \eqref{def:Perturb:F}, \eqref{def:Perturb:G} and \eqref{def:Perturb:R}. Thanks to Hypothesis \textbf{HP2}  the functions $F$ and $R$ are trigonometric polynomials of degree at most $M$. Then $F$, which by construction has zero average, can be written as
\begin{equation*}
 F(x,\tau)=\sum_{k=-M}^M a_k(\tau)e^{ikx}
\end{equation*}
for certain analytic periodic functions $a_k:\TT\rightarrow \CC$ with zero average. We write the function $R$ in \eqref{def:Perturb:R} as
\begin{equation*}
 R(x,\eps y,\tau;\eps)=\sum_{\substack{k=-M\ldots M\\ l\geq 1}}c_{kl}(\tau;\eps)e^{ikx}(\eps y)^l
\end{equation*}
for certain analytic functions $c_{kl}:\TT\rightarrow \CC$. Finally, we write $G$ in \eqref{def:Perturb:G} as
\begin{equation*}
G(\eps y)=\sum_{k\geq 3}g_k (\eps y)^k
\end{equation*}
for certain constants $g_k\in\RR$.

By Hypothesis \textbf{HP1},  the
Hamiltonian $H$ in \eqref{def:ToyModel} is analytic in $\tau=t/\eps$.
By the compactness of $\TT$, actually  there exists a constant
$\sigma_0$ such that $H$ is continuous in $\ol \TT_{\sigma_0}$ and
analytic in $\TT_{\sigma_0}$.  From now on, we fix  $0<\sigma<\sigma_0$.

Throughout the proof of  Theorem
\ref{th:Main} we will use the analyticity in $\mu$. We fix an
arbitrary value  $\mu_0>0$. Even if we do not write it explicitly, all
functions which appear in this paper will be analytic in $\mu \in B(\mu_0)=\{\mu\in \CC:|\mu|<\mu_0\}$.

From now on, we work with the fast time $\tau=t/\eps$. Then,
denoting $'=d/d\tau$, we have the system
\begin{equation}\label{eq:ode:original:lent}
\left\{\begin{array}{rl} x'&\dps=\eps\left(y+\mu \pa_y
H_1\left(x,y,\tau;\eps\right)\right)\\
y'&\dps =-\eps\left(V'(x)+\mu\pa_x
H_1\left(x,y,\tau;\eps\right)\right).
\end{array}\right.
\end{equation}


\subsection{The periodic orbit}
The next theorem
states the existence and useful properties of a hyperbolic periodic
orbit close to the origin of the perturbed system. 

\begin{theorem}\label{th:Periodica}
Assume Hypotheses \textbf{HP1} and  \textbf{HP4}  and  fix any $\mu_0>0$. Then,
there exists $\eps_0>0$ such that for any $|\mu|<\mu_0$ and
$\eps\in(0,\eps_0)$,  system \eqref{eq:ode:original:lent} has a
$2\pi$-periodic orbit $(\xp(\tau),\yp(\tau)):\TT_\sigma\rightarrow
\CC^2$ which is real-analytic and satisfies that for $\tau\in \TT_\sigma$,
\[
\begin{split}
\left|\xp(\tau)\right|&\leq
b_0|\mu|\eps^{2}\\
\left|\yp(\tau)\right|&\leq
b_0|\mu|\eps,
\end{split}
\]
where $b_0>0$ is a constant independent of $\eps$ and $\mu$.
\end{theorem}

The proof of this theorem follows the same lines as the proof of Theorem 4.1 in \cite{BaldomaFGS11}. Note that this theorem gives the same result as Theorem \ref{th:MainPO}. The only difference is that it proves the existence of the periodic orbit in an analytic strip around the real axis. This is necessary to deal later on with analytic parameterizations of the invariant manifolds.

\begin{remark}
The periodic orbit $(\xp(\tau),\yp(\tau))$, 
 the Hamiltonians $\wh H$ and 
$\wh H_1$, which will be defined in \eqref{def:HamPeriodicaShiftada} and \eqref{def:ham:ShiftedOP:perturb} respectively, and the functions involved in these definitions depend on the
parameters $\mu$, $\eps$.
From now on, we omit this dependence explicitly.
\end{remark}

Once we know the existence of the periodic orbit, we perform the
time dependent change of variables
\begin{equation}\label{eq:CanviShiftPO}
\left\{\begin{array}{l} q=x-\xp(\tau)\\ p=y-\yp(\tau)
\end{array}\right.
\end{equation}
which transforms system \eqref{eq:ode:original:lent} into a
Hamiltonian system with Hamiltonian function $\eps\widehat
H(q,p,\tau)$:
\begin{equation}\label{def:HamPeriodicaShiftada}
\begin{split}
\wh
H(q,p,\tau)=&\frac{p^2}{2}
+V\left(q+\xp(\tau)\right)-V\left(\xp(\tau)\right)-V'\left(\xp(\tau)\right)q\\
&+\mu\widehat H_1(q,p,\tau)
\end{split}
\end{equation}
where
\begin{equation}\label{def:ham:ShiftedOP:perturb}
 \wh H_1(q,p,\tau)=\eps^{-2}\wh G(\eps p,\tau)+\wh F(q,\tau)+\wh R(q,\eps p,\tau)
\end{equation}
with
\begin{align*}
\wh F(q,\tau)= &F(q+x_p(\tau),\tau)-F(x_p(\tau),\tau)-\pa_x F(x_p(\tau),\tau)q\\
\wh R(q,\eps p,\tau)=& R(q+x_p(\tau),\eps p+\eps y_p(\tau),\tau)-R(x_p(\tau),\eps y_p(\tau),\tau)\\
&-\pa_x R(x_p(\tau),\eps y_p(\tau),\tau)q-\pa_y R(x_p(\tau),\eps y_p(\tau),\tau)\eps p\\
\wh G(\eps p,\tau)=&G\left(\eps p+\eps y_p(\tau)\right)-G\left(\eps y_p(\tau)\right)-\eps G'\left(\eps y_p(\tau)\right)p
\end{align*}
We have added the terms
$V\left(\xp(\tau)\right)$, $F(x_p(\tau),\tau)$, $R(x_p(\tau),\eps y_p(\tau),\tau)$ and $G\left(\eps y_p(\tau)\right)$ for convenience.
Note that they do not generate any term in the differential equations associated
to $\wh H$. Since
$|x_p(\tau)|=\OO\left(\mu \eps^2\right)$, we can split the function $\wh F$ as 
\[
\wh F(q,\tau)=\wh F_1(q,\tau)+ \wh F_2(q,\tau),
\]
where
\[
\wh F_1(q,\tau)=\sum_{k=-M\ldots M}a_k(\tau) \left(e^{kiq}-1-kiq\right)
\]
and $ \wh F_2(q,\tau)$ is the remaining part, which satisfies $\wh F_2(q,\tau)=\OO(\mu\eps^2)$. Recall that the functions $a_k(\tau)$ have zero average and therefore we have that 
\[
 \langle F_1\rangle=0.
\]

\subsection{Existence of parameterizations of the perturbed invariant manifolds}\label{sec:DiferentsParam}
The next step is to prove the existence of parameterizations of the unstable and stable invariant
manifolds of the periodic orbit given in Theorem \ref{th:Periodica}.  

As is done in \cite{BaldomaFGS11} (see also \cite{Baldoma06}), we follow \cite{LochakMS03,
Sauzin01}, and we write the invariant manifolds as graphs of suitable generating
functions which are solutions of a Hamilton-Jacobi equation in appropriate variables.
To this end, we consider the
symplectic change of variables (see \cite{Baldoma06, BaldomaFGS11})
\begin{equation}\label{eq:CanviSimplecticSeparatriu}
\left\{\begin{array}{l}
q=q_0(u)\\
\dps p=\frac{w}{p_0(u)},
\end{array}\right.
\end{equation}
where $(q_0(u),p_0(u))$ is the parameterization of the homoclinic orbit given in
Hypothesis \textbf{HP6}. This change is well defined for any $u\in\CC$ such that $p_0(u)\neq 0$ and
leads to a new Hamiltonian given  by
\begin{equation}\label{def:Hbarra}
\eps \overline H(u,w,\tau)=\eps \widehat
H\left(q_0(u),\frac{w}{p_0(u)}, \tau\right),
\end{equation}
where $\widehat H$ is the Hamiltonian defined in
\eqref{def:HamPeriodicaShiftada}. When $\mu=0$, $\wh H$ becomes $H_0$ defined in
\eqref{def:ToyModel:H0}. Then, the separatrix of the unperturbed
system ($\mu=0$) for $\ol H$ can be parameterized as a graph as
$w=p_0(u)^2$. If we want to obtain parameterizations of the perturbed invariant
manifolds, we can take into account the well known fact that,  locally,  they
are Lagrangian and can be obtained as graphs of generating functions which are
solutions of the Hamilton-Jacobi equation  associated to the Hamiltonian $\eps\overline H$.
That is, we parameterize the invariant manifolds as graphs as $w= \partial _{u} T^{u,s}(u,\tau)$, where the
functions $T^{u,s}$ satisfy
\begin{equation}\label{eq:HamJacGeneral}
\pa_\tau T(u,\tau)+ \eps \overline H(u,\pa_u T(u,\tau),\tau)=0.
\end{equation}
The solutions of this equation give parameterizations of the
invariant manifolds, which, in the original variables, read
\begin{equation}\label{eq:ParameterizationHJ}
(q,p)=\left(q_0(u), \frac{\pa_u T^{u,s}(u,\tau)}{p_0(u)}\right).
\end{equation}
The inconvenient of these parameterizations is that they are not defined everywhere since $p_0(u)$ might vanish in some  $u\in \CC$. In some concrete cases, as happens for the classical pendulum, $p_0(u)\neq 0$ for any $u\in\CC$, and therefore we do not have any problem. However, this is not always de case. In particular, for the figure eight separatrices (see Figure \ref{fig:FigureEight}) there is always a real value of $u$ such that $p_0(u)=0$. For the graph case (see Figure \ref{fig:Pendol}), $p_0(u)\neq 0$ for all $u\in\RR$ but it can vanish for complex values of the variable. Therefore, as was done in \cite{BaldomaFGS11}, we look for parameterizations of the form \eqref{eq:ParameterizationHJ} in very special complex domains for the variable $u$, which were defined in that paper and were called \emph{boomerang domains} (see Figure \ref{fig:BoomerangDomains}). The reason of the somewhat strange shape of these domains is explained in full detail in that paper. Let us just say here that we need domains which contain both an interval of the real line of width independent of $\eps$ and points at a distance $\OO(\eps)$ of the singularities $u=\pm ia$. Moreover, we want these properties to be still satisfied when we intersect the boomerang domains associated to the stable and unstable manifolds.  These domains are defined as follows
\begin{equation}\label{def:DominisRaros}
\begin{split}
D^{s}_{\kk,d}=&\left\{u\in\CC;\right.\left. |\Im u|<\tan\beta_1\Re
u+a-\kk\eps,
|\Im u|<\tan\beta_2\Re u+a-\kk\eps,\right.\\
& \left.|\Im u|>\tan\beta_2\Re u+a-d\right\}\\
D^u_{\kk,d}=&\left\{u\in\CC;\right.\left.|\Im u|<-\tan\beta_1\Re
u+a-\kk\eps,
|\Im u|<\tan\beta_2\Re u+a-\kk\eps,\right.\\
& \left.|\Im u|>\tan\beta_2\Re u+a-d\right\}\\
&\cup \left\{ u\in\CC; \right.|\Im u|<-\tan\beta_1\Re
u+a-\kk\eps, |\Im u|>-\tan\beta_2\Re u+a-d,\\
& \left.\Re u<0\right\},
\end{split}
\end{equation}
where $\beta_1\in (0,\pi/2)$ is any fixed angle. The choice of $\beta_2$ goes as follows.
First let us point out that the zeros of $p_0(u)$ are isolated in
$\CC$. Moreover, close to the singularities $u=\pm ia$, $p_0(u)$ can
not vanish. Then, in order to assure that $p_0(u)$ does not vanish
in the whole domains $D^s_{\kk,d}$ and $D^u_{\kk,d}$, one has to
choose an angle $\beta_2$ such that $\beta_2>\beta_1$ has a positive
lower bound independent of $\eps$ and $\mu$  and such that the lines $|\Im
u|=\tan\beta_2\Re u+a$ do not contain any zero of $p_0(u)$. Then,
taking $\eps>0$ and $d>0$ independent of $\eps$, both small enough,
one can guarantee that $p_0(u)$ does not vanish neither in
$D^s_{\kk,d}$ nor in $D^u_{\kk,d}$.

\begin{figure}[H]
\begin{center}
\psfrag{u1}{$u_1$}\psfrag{u2}{$\bar u_1$}
\psfrag{b1}{$\beta_1$}\psfrag{b2}{$\beta_2$}\psfrag{a1}{$ia$}
\psfrag{a2}{$-ia$}\psfrag{a3}{$i(a-d)$}\psfrag{a4}{$i(a-\kk\eps)$}
\psfrag{D1}{$D^{\out,s}_{\rr,
\kk}$}\psfrag{D}{$D^{s}_{\kk,d}$}\psfrag{D4}{$D^{\out,u}_{\rr,\kk}$}\psfrag{D3}{
$D^{u}_{\kk,d}$}
\includegraphics[height=6cm]{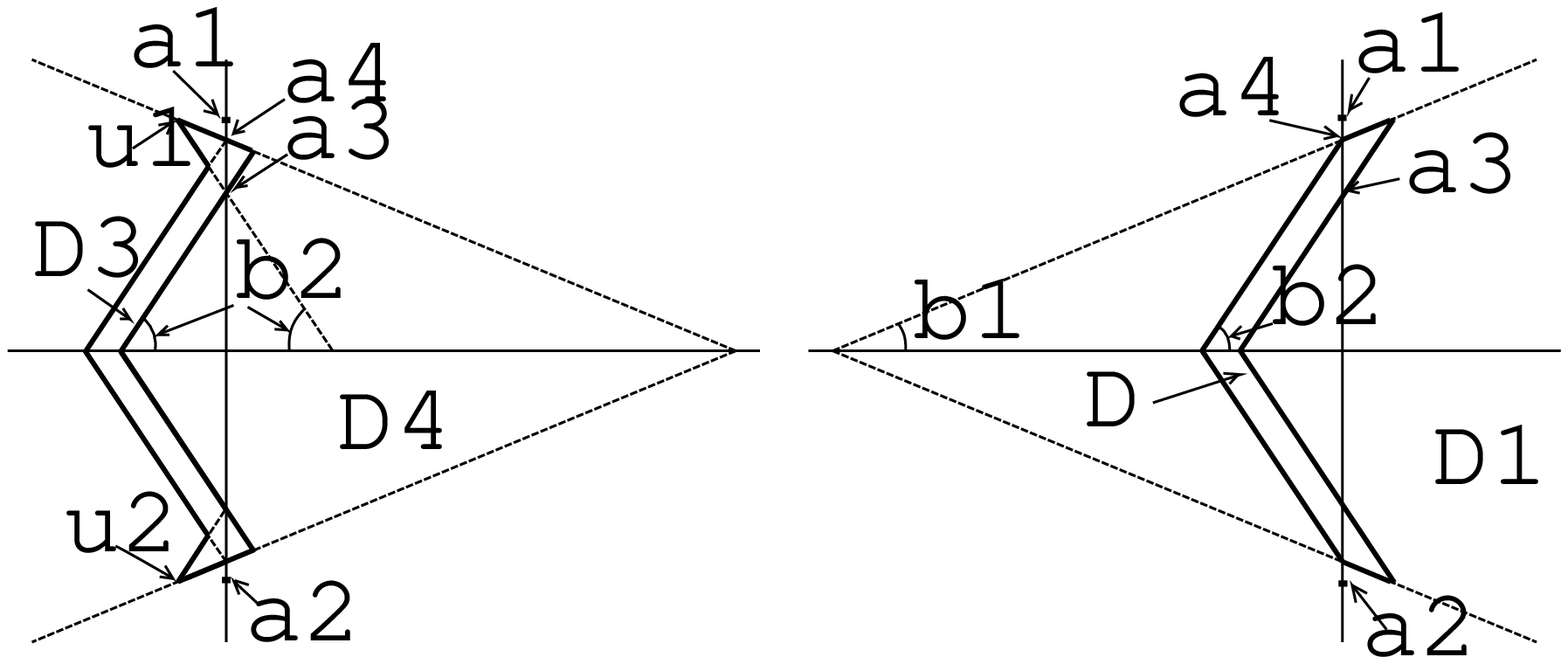}
\end{center}
\caption{\figlabel{fig:BoomerangDomains} The \emph{boomerang
domains} $D^{u}_{\kk,d}$ and $D^{s}_{\kk,d}$ defined in
\eqref{def:DominisRaros}.}
\end{figure}
As it is well known, to obtain the parameterizations of the invariant manifolds using  equation \eqref{eq:HamJacGeneral}, one has to impose the asymptotic conditions 
\begin{align}
\dps \lim_{\Re u\rightarrow-\infty}p_0\ii(u)\cdot\pa_u T^u(u,\tau)=0 & \text{\quad (for the unstable manifold)} \label{eq:AsymptCondFuncioGeneradora:uns}\\
\dps \lim_{\Re u\rightarrow+\infty}p_0\ii(u)\cdot \pa_u
T^s(u,\tau)=0 & \text{\quad (for the stable manifold)}.
\label{eq:AsymptCondFuncioGeneradora:st}
\end{align}
These asymptotic conditions
ensure that the invariant manifolds tend to the periodic orbit 
$(q,p)=(0,0)$ of the system associated to Hamiltonian \eqref{def:HamPeriodicaShiftada} as $\Re u\rightarrow\pm\infty$. 
Of course these conditions do not have any meaning in the domains $D^{u,s}_{\kk,d}$ since these domains are bounded. This implies that, to prove the existence of the parameterizations of the invariant manifolds in these domains, one has to start with different domains were these asymptotic conditions make sense and then one has to find a way to extend them analytically to the domains $D^{u,s}_{\kk,d}$ taking into account that these parameterizations become undefined at the points where $p_0(u)=0$. This process is explained in great detail in \cite{BaldomaFGS11}. Following the notation of that paper, the case we are considering now ressembles the case $r=1$ and $\ell=2$ of that paper. Therefore here we just sketch it and give the main ideas. 

The main steps of this procedure are the following.
\begin{itemize}
 \item[1] We first obtain parameterizations of the form \eqref{eq:ParameterizationHJ} satisfying the asymptotic conditions \eqref{eq:AsymptCondFuncioGeneradora:uns} and \eqref{eq:AsymptCondFuncioGeneradora:st} in the domains 
\begin{equation*}
\begin{array}{l}
D^{u}_{\infty,\rr}=\{u\in\CC; \Re u<-\rr\}\\
D^{s}_{\infty,\rr}=\{u\in\CC; \Re u>\rr\}.
\end{array}
\end{equation*}
Note that thanks to the exponential decay of $p_0(u)$ as $\Re u\rightarrow\pm\infty$, one can ensure that for $\rr$ big enough, $p_0(u)\neq 0$ for $u\in D^{u,s}_{\infty,\rr}$.
\item[2] To reach the boomerangs domains $D^{u,s}_{\kk,d}$, we cannot extend the just obtained parameterizations since in the analytical extension procedure, we might reach points where $p_0(u)=0$.  Thus, we switch to a different kind of parameterizations of the invariant manifolds, which are chosen of the form 
\begin{equation}\label{def:ParamByFlow}
(q,p)=(Q(v,\tau),P(v,\tau))
\end{equation}
in such a way that $(Q(v+\eps s,\tau+s),P(v+\eps s,\tau+s))$
are solutions of the differential equation associated to  Hamiltonian
\eqref{def:HamPeriodicaShiftada}. We will extend these parameterizations and once we reach the  domains $D^{u,s}_{\kk,d}$, in which $p_0(u)\neq0$, we will switch back to parameterizations of the form \eqref{eq:ParameterizationHJ}. This means that we need to define the parameterizations \eqref{def:ParamByFlow} in domains which overlap both with the domains $D^{u,s}_{\infty,\rr}$ and with the domains $D^{u,d}_{\kk,d}$. We define the following domains (see also Figure \ref{fig:OuterParam}),
\begin{equation}\label{def:DominOuterParam}
\begin{split}
D^{\out,u}_{\rr',d,\kk}=\Big\{u\in\CC
; \,\,&|\Im u|<-\tan \beta_1\Re u+a-\kk\eps,\, \Re u>-\rr',\\
&|\Im
u|<-\tan\beta_2\Re u+a-\frac{d}{2}\Big\}\\
D^{\out,s}_{\rr',d,\kk}=\Big\{u\in\CC;\,\,& |\Im u|<\tan \beta_1\Re u+a-\kk\eps,\, \Re u<\rr'\\
& |\Im u|<\tan\beta_2\Re u+a-\frac{d}{2}\Big\},
\end{split}
\end{equation}
where $\rr'>\rr$ so that $D^{\out,\ast}_{\rr,d,\kk}\cap D^\ast_{\infty,\rr}\neq\emptyset$, $\ast=u,s$.

\begin{figure}[h]
\begin{center}
\psfrag{D6}{$D^{\out,s}_{\rr,d,\kk}$} \psfrag{D5}{$
D^{\out,u}_{\rr,d,\kk}$}
\psfrag{b1}{$\beta_1$}\psfrag{b2}{$\beta_2$}\psfrag{r1}{$-\rr$}\psfrag{r2}{$\rr$
}
\psfrag{a1}{$ia$}\psfrag{a2}{$-ia$}\psfrag{a3}{$i(a-d)$}\psfrag{a4}{
$i(a-\kk\eps)$}
\psfrag{D1}{}\psfrag{D}{$D^{s}_{\kk,d}$}\psfrag{D4}{}\psfrag{D3}{$D^{u}_{\kk,d}$}
\includegraphics[height=6cm]{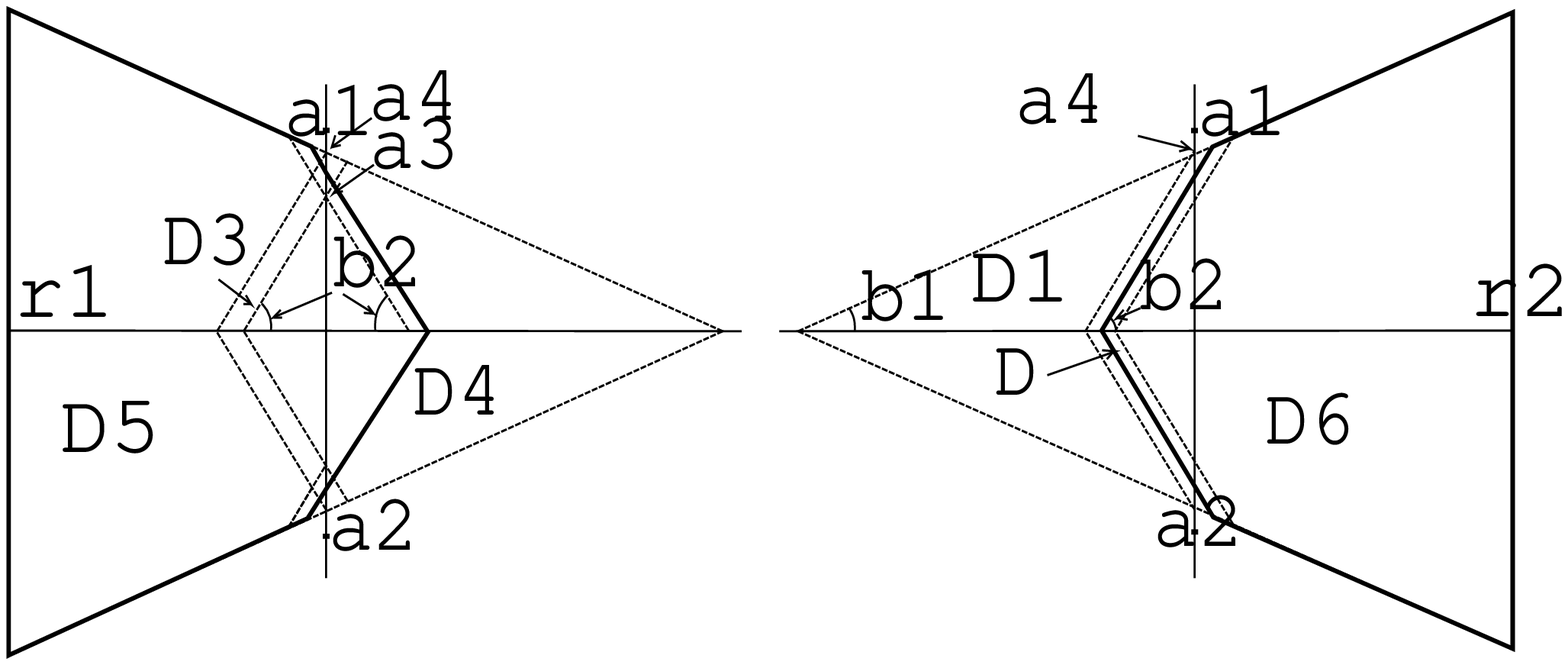}
\end{center}
\caption{\figlabel{fig:OuterParam} The domains $
D^{\out,u}_{\rr,d,\kk}$ and  $D^{\out,s}_{\rr,d,\kk}$ defined in
\eqref{def:DominOuterParam}.}
\end{figure}
The parameterizations \eqref{def:ParamByFlow} were used in
\cite{DelshamsS92, DelshamsS97,Gelfreich97a,Gelfreich00,BaldomaF04,BaldomaF05}.
It is easy to see
\cite{Gelfreich97a} that $(Q,P)$ have to satisfy
\begin{equation}\label{eq:PDEParametritzacions}
\LL_\eps\left(\begin{array}{c}Q\\P\end{array}\right)=\left(\begin{array}{c}
P+\mu\pa_p \widehat H_1(Q,P,\tau)\\
-\left(V'(Q+\xp(\tau))-V'(\xp(\tau))\right)-\mu\pa_q
\widehat H_1(Q,P,\tau)\end{array}\right),
\end{equation}
where $\LL_\eps$ is the operator
\begin{equation*}
\LL_\eps=\eps\ii\pa_\tau+\pa_v
\end{equation*}
and $\widehat H_1$ is the Hamiltonian defined in
\eqref{def:ham:ShiftedOP:perturb}.
\item[3] To make this analytic extension procedure with parameterizations of the type \eqref{def:ParamByFlow}, we first derive parameterizations $(Q^{u,s}(v,\tau),P^{u,s}(v,\tau))$ in  the  domains $D^{\out,u}_{\kk,\rr'}\cap
D^{u}_{\infty,\rr}$ and $D^{\out,s}_{\kk,\rr'}\cap
D^{s}_{\infty,\rr}$ from the parameterizations \eqref{eq:ParameterizationHJ}  obtained in Step 1. Taking into
account the change of variables
\eqref{eq:CanviSimplecticSeparatriu}, it is natural to look for these
parameterizations as
\begin{equation*}
\begin{array}{ll}
\dps Q^{u,s}(v,\tau)=q_0\left(v+\UU^{u,s}(v,\tau)\right)\\
\dps P^{u,s}(v,\tau)=\frac{\pa_u
T^{u,s}\left(v+\UU^{u,s}(v,\tau)\right)}{p_0(v+\UU^{u,s}(v,\tau))},
\end{array}
\end{equation*}
where $\UU^{u,s}$ define changes of variables
$u=v+\UU^{u,s}(v,\tau)$  so that $(Q^{u,s},P^{u,s})$
satisfy the system of equations \eqref{eq:PDEParametritzacions}.
\item[4] Once we have obtained these parameterizations, we  use equation \eqref{eq:PDEParametritzacions} to extend them to the whole domains $D^{\out,u}_{\rr,d,\kk}$ and $D^{\out,s}_{\rr,d,\kk}$. This step is straightforward since these domains are not close to the singularities $u=\pm ia$. 
\item[5] Now we  have  parameterizations of the
invariant manifolds of the form \eqref{def:ParamByFlow} in the
following  domains, which are depicted  in Figure \ref{fig:TransBoom},
\begin{equation}\label{def:Dominis:Trans:Raros}
\begin{split}
\tro^{\out,u}_{\kk,d}&=D^{\out,u}_{\rr,d,\kk}\cap D^{u}_{\kk,d}\\
\tro^{\out,s}_{\kk,d}&=D^{\out,s}_{\rr,d,\kk}\cap D^{s}_{\kk,d},
\end{split}
\end{equation}
where, by construction, $p_0(u)$ does not vanish. Thus, we can use
these domains as  transition domains where we can go back to the
parameterizations \eqref{eq:ParameterizationHJ} and where the
Hamilton-Jacobi equation \eqref{eq:HamJacGeneral} can be used. To
obtain them, we look for changes of variables
$v=u+\VV^{u,s}(u,\tau)$ which satisfy
\begin{equation*}
Q^{u,s}(u+\VV^{u,s}(u,\tau),\tau)=q_0(u),
\end{equation*}
where $Q^{u,s}$ are the first components of the parameterizations
just obtained. From these changes, we can deduce the generating functions $T^{u,s}$ which  give the
parameterizations \eqref{eq:ParameterizationHJ}. 
\begin{figure}[h]
\begin{center}
\psfrag{D6}{$\wt D^{\out,s}_{\rr,d,\kk}$} \psfrag{D5}{$\wt
D^{\out,u}_{\rr,d,\kk}$}
\psfrag{b1}{$\beta_1$}\psfrag{b2}{$\beta_2$}
\psfrag{a1}{$ia$}\psfrag{a2}{$-ia$}\psfrag{a3}{$i(a-d)$}\psfrag{a4}{
$i(a-\kk\eps)$}
\psfrag{D1}{$D^{\out,s}_{\rr,\kk}$}\psfrag{D}{$D^{s}_{\kk,d}$}\psfrag{D4}{$D^{
\out,u}_{\rr,\kk}$}\psfrag{D3}{$D^{u}_{\kk,d}$}
\psfrag{D8}{$\tro^{\out,s}_{\kk,d}$}\psfrag{D7}{$\tro^{\out,u}_{\kk,d}$}
\includegraphics[height=6cm]{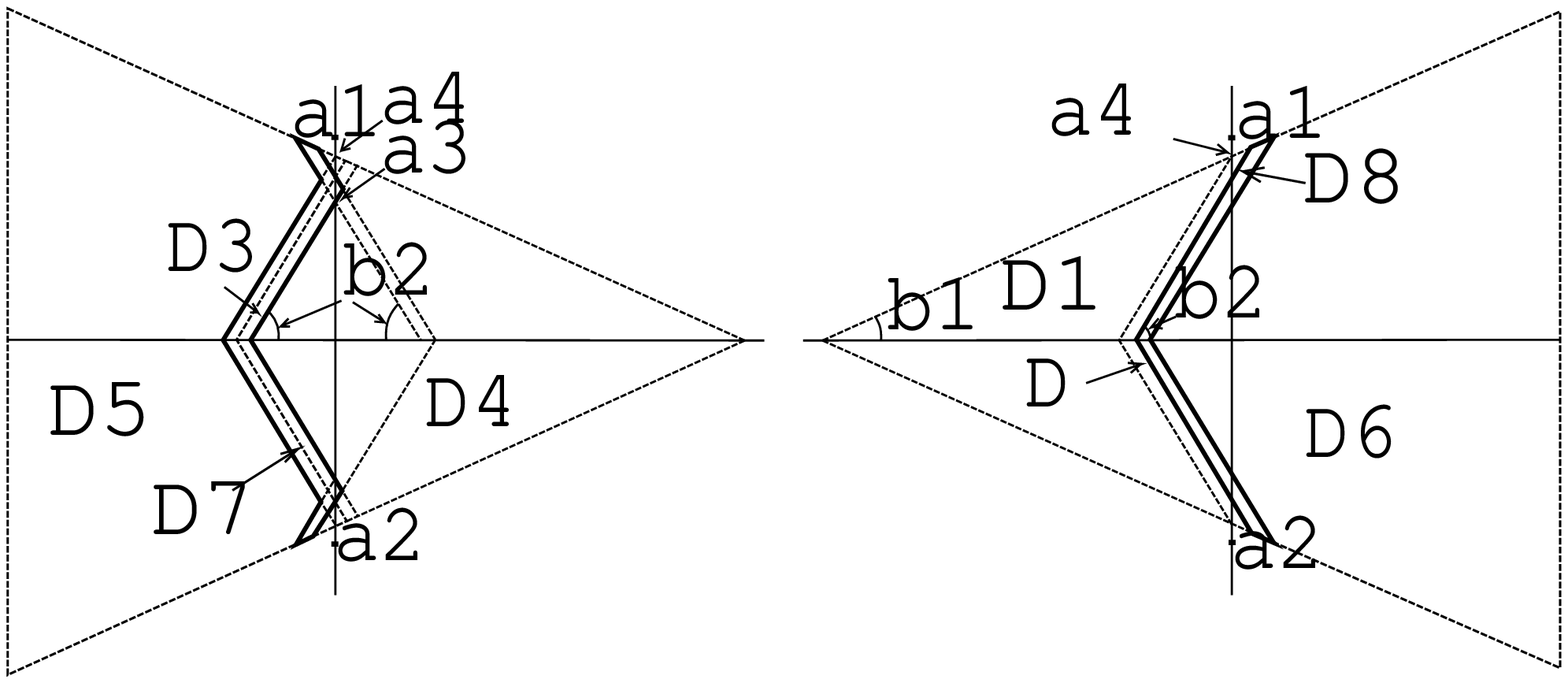}
\end{center}
\caption{\figlabel{fig:TransBoom} The domains
$\tro^{\out,u}_{\kk,d}$ and $\tro^{\out,s}_{\kk,d}$ defined in
\eqref{def:Dominis:Trans:Raros}.}
\end{figure}
\item[6] The last step is to extend analytically the parameterizations \eqref{eq:ParameterizationHJ} just obtained   to the whole domains  $D^{u,s}_{\kk,d}$. Since in these domains $p_0(u)\neq0$, we take into account that the associated generating functions   $T^{u,s}$ satisfy equation \eqref{eq:HamJacGeneral} and we use this equation to obtain these final analytic extensions.
\end{itemize}
Note that in the case that $p_0(u)\neq 0$ for any $u\in\CC$, as happens for the classical pendulum, one can skip steps 2 to 5. That is, one can extend the parameterizations \eqref{eq:ParameterizationHJ} from infinity directly to the boomerang domains $D^{u,s}_{\kk,d}$.

As a conclusion of this procedure, which is explained in fully detail in Sections 6 and 7 of \cite{BaldomaFGS11}, we can state the following theorem. Recall that using the notation of that paper, the problem considered in the present paper is stronly related to the case $r=1$ and $\ell=2$ considered in that paper. We state the result for the unstable invariant manifold. The stable one satisfies analogous properties.
\begin{theorem}\label{th:ExtensioFinal}
Let $\mu_0$  be the constant given in Theorem
\ref{th:Periodica},  $d_1>0$, $\kk_1>0$ big enough and
$\eps_0>0$ small enough. Then, for $\mu\in B(\mu_0)$ and
$\eps\in(0,\eps_0)$, the unstable invariant manifold of the periodic orbit $(q,p)=0$ of the system associated to Hamiltonian \eqref{def:HamPeriodicaShiftada} has a parameterization of the form  \eqref{eq:ParameterizationHJ} for $(u,\tau)\in D^{u}_{\kk_1,d_1}\times\TT_\sigma$.

Moreover, there exists a real constant $b_1>0$ independent of $\eps$
and $\mu$ such that, for $(u,\tau)\in
D^{u}_{\kk_1,d_1}\times\TT_\sigma$, the generating function $T^u$ which gives the parameterization \eqref{eq:ParameterizationHJ} satisfies that
\[
\left|\pa_uT^u(u,\tau)-\pa_u T_0(u)\right|\leq
\frac{b_1|\mu|\eps}{\left|u^2+a^2\right|^{3}},
\]
where 
\begin{equation}\label{def:T00}
T_0(u)=\int_{-\infty}^u p_0^2(v)\,dv,
\end{equation}
is the parameterization of the unperturbed separatrix.
\end{theorem}


\subsection{The asymptotic first order of $\partial_u T^{u,s}$ close to the
singularities $\pm i a$}\label{sec:Aproximacio}

Theorem \ref{th:ExtensioFinal} gives parameterizations of the invariant manifolds up to points at a distance of order $\OO(\eps)$ of the singularities. Nevertheless, at this distance of the singularities Theorem \ref{th:ExtensioFinal} is not a perturbative result with respect to the singular parameter $\eps$. Namely, at a distance $\OO(\eps)$ of $u=\pm ia$  the parameterizations of the perturbed invariant manifolds are not well approximated by the parameterization of the unperturbed separatrix. This fact can be easily seen from \eqref{eq:SepartriuAlPolTrig} and the bounds given in Theorem \ref{th:ExtensioFinal}, since for $u\mp ia\sim\eps$ both  $\pa_u T_0(u)=p_0^2(u)$ and the remainder become of size $\OO(1/\eps^2)$. Then, to study the difference between the manifolds, we need to look for better approximations of $T^{u,s}$ in the following domains, which are usually called \emph{inner domains}  (see Figure
\ref{fig:Inners}),
\begin{equation}\label{def:DominisInnerEnu}
\begin{split}
D_{\kk,\C}^{\inn,+,u}=&\left\{ u\in\CC; \Im
u>-\tan\beta_1 (\Re u+\C\eps^\gamma)+a, \Im u<-\tan \beta_2\Re
u+a-\kk\eps,\right.\\
& \left.\Im u<-\tan\beta_0\Re u+a-\kk\eps\right\}\\
D_{\kk,\C}^{\inn,-,u}=&\left\{u\in\CC; \bar u \in
D_{\kk,\C}^{\inn,+,u}\right\}\\
D_{\kk,\C}^{\inn,+,s}=&\left\{u\in\CC; -\bar u \in
D_{\kk,\C}^{\inn,+,u}\right\}\\
D_{\kk,\C}^{\inn,-,s}=&\left\{u\in\CC; - u \in
D_{\kk,\C}^{\inn,+,u}\right\}
\end{split}
\end{equation}
where $\kk>0$, $\C>0$ and $\gamma \in (0,1)$. The angles $\beta_1$  and $\beta _2$ are the ones  considered in
the definition of the boomerang domains in \eqref{def:DominisRaros}  and
$\beta_0$ is  any angle
satisfying that $\beta_1-\beta_0$ has a positive lower bound independent of
$\eps$ and $\mu$. Let us observe that, if $u\in D_{\kk,\C}^{\inn,\pm,\ast}$, $\ast=u,s$, and $\eps$ is small enough, then 
$\OO(\kk \eps) \leq |u\mp ia| \leq \OO(\eps^{\gamma})$ and  $u\in  D^\ast_{\kk,d}$,  $\ast=u,s$. 

\begin{figure}[H]
\begin{center}
\psfrag{D}{$D_{\kk,\C}^{\inn,+,u}$}\psfrag{D1}{$D_{\kk,\C}^{\inn,-,u}$}
\psfrag{D3}{$D_{\kk,\C}^{\inn,+,s}$}\psfrag{D4}{$D_{\kk,\C}^{\inn,-,s}$}
\psfrag{b0}{$\beta_0$}\psfrag{b1}{$\beta_1$}\psfrag{b2}{$\beta_2$}
\psfrag{a}{$ia$}\psfrag{a1}{$-ia$}\psfrag{a2}{$i(a-\kk\eps)$}
\includegraphics[height=6cm]{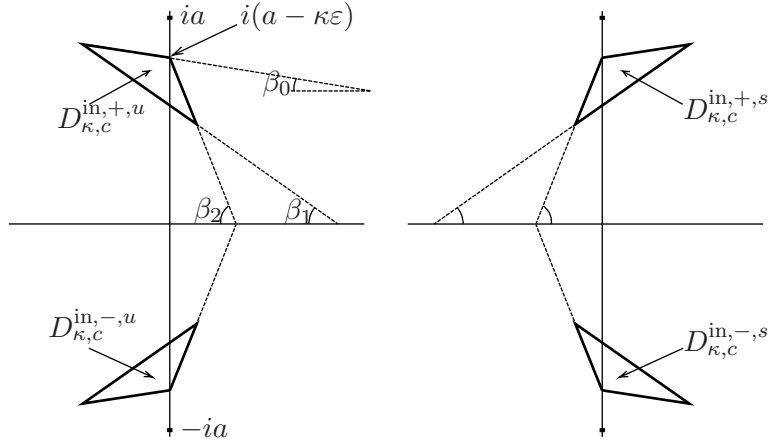}
\end{center}
\caption{\figlabel{fig:Inners} The \emph{inner domains} defined in
\eqref{def:DominisInnerEnu}.}
\end{figure}

We obtain good approximations of the functions $T^{u,s}$ in these domains  through a singular
limit.  Since the study of both invariant
manifolds close either to $u=ia$ or $u=-ia$ is analogous, we only
study them in the domain $D_{\kk,\C}^{\inn,+,u}$. Thus, we consider
the change of variables
\begin{equation*}
z=\eps\ii (u-ia).
\end{equation*}
The variable $z$ is called the \emph{inner variable}, in
contraposition to the \emph{outer variable} $u$. By
definition of $T_0$ in \eqref{def:T00} and using the expansion
around the singularities of $p_0(u)$ in \eqref{eq:SepartriuAlPolTrig}, we have that
\begin{equation*}
\partial_uT_0(\eps z+ ia) =p_0^2\left(\eps z+ia\right)=\frac{C_+^{2}}{\eps^{2}z^{2}}
\left(1+\OO\left((\eps z)^{2/M}\right)\right)
\end{equation*}
and, using the results of Theorem \ref{th:ExtensioFinal}, we have that
\begin{equation*}
\left|\partial_u T^{u,s}(\eps z+ ia, \tau) - \partial_u T_0(\eps z+
ia)\right|\leq K\frac{|\mu|}{\eps^{2}|z|^{3}}.
\end{equation*}
Hence, to study the first order in $\eps$ of the functions $T^{u,s}$ at a distance $\OO(\eps)$ of the singularities we
scale the generating function as
\begin{equation}\label{eq:FuncioGeneradoraInner}
\psi^{u,s}(z,\tau)=\eps C_+^{-2} T^{u,s}(ia+\eps z,\tau).
\end{equation}
Then, the Hamilton-Jacobi equation \eqref{eq:HamJacGeneral} becomes
\begin{equation}\label{eq:HJGeneralInner}
\pa_\tau\psi +\eps^{2}C_+^{-2}\ol H\left(ia+\eps z,
\eps^{-2}C^{2}_+\pa_z\psi,\tau\right)=0,
\end{equation}
where $\ol H$ is the Hamiltonian function defined in
\eqref{def:Hbarra}. The corresponding Hamiltonian is
\begin{equation}\label{Hamiltonia:varInner}
\HH (z,w,\tau)=\eps^{2}C_+^{-2}\overline H\left(ia+\eps z,
\eps^{-2}C^{2}_+w,\tau\right).
\end{equation}
We study equation \eqref{eq:HJGeneralInner} in the domain
$\DD_{\kk,\C}^{\inn,+,u}\times\TT_\sigma$, where
\begin{equation}\label{def:DominisInnerEnz}
\begin{split}
\DD_{\kk,\C}^{\inn,+,u}=&\left\{ z\in\CC; ia+\eps z \in
D_{\kk,\C}^{\inn,+,u}\right\}.
\end{split}
\end{equation}

To study equation \eqref{eq:HJGeneralInner}, as a first step it is
natural to study it in the limit case $\eps=0$, which reads
\begin{equation}\label{eq:HJEqInner}
\pa_\tau\psi_0
+\frac{1}{2}z^{2}\left(\pa_z\psi_0\right)^2-\frac{1}{2z^{2}}+\mu\left(\wt G\left(z\pa_z\psi_0\right)+
\frac{\wt a(\tau)}{z^2}+\frac{1}{z^2}\sum_{l\geq 1}\wt c_l(\tau)\left(z\pa_z\psi_0\right)^l\right)=0,
\end{equation}
where
\[
\wt G(w)=C_+^{-2} G(C_+w)
\]
The other terms are the ones coming from the functions $F$ and $R$ defined in \eqref{def:Perturb:F} and \eqref{def:Perturb:R} respectively. Indeed,  $\mu\wt a(\tau)/z^2$ is the leading term coming from the degree $M$  terms of $F(x,\tau)$  (recall that $F$ is a trigonometric polynomial in $x$ of degree $M$ thanks to  Hypothesis \textbf{HP2}),  and  
\[
\mu \frac{1}{z^2}\sum_{l\geq 1}\wt c_l(\tau)\left(z\pa_z\psi_0\right)^l
\]
comes from the degree $M$ terms of the function $R(x,I,\tau)$ in \eqref{def:Perturb:R}.

Equation \eqref{eq:HJEqInner} is the Hamilton-Jacobi equation associated to the non-autonomous Hamiltonian
\[
\HH_0(z,w,\tau)=\frac{1}{2}z^{2}w^2-\frac{1}{2z^{2}}+\mu\left(\wt G(zw)+
\frac{\wt a(\tau)}{z^2}+\frac{1}{z^2}\sum_{l\geq 1}\wt c_l(\tau)\left(zw\right)^l\right),
\]
which satisfies that $\HH\rightarrow\HH_0$ as $\eps\rightarrow 0$,
where $\HH$ is the Hamiltonian function defined in
\eqref{Hamiltonia:varInner}.

\begin{remark}\label{remark:Example:Inner}
As an example we show the inner equation associated to the Hamiltonian \eqref{def:Example:Ham:rescaled} considered in Section \ref{sec:Example}. Recall that, as we have seen in Section  \ref{sec:Example}, the integrable system $H_0$ associated to \eqref{def:Example:Ham:rescaled} is simply the classical pendulum, whose upper separatrix has singularities at $u=\pm i\pi/2 ki\pi$, $k\in \ZZ$ and $C_+=-2i$. Then, the inner equation around $u=i\pi/2$ is given by
\[
 \pa_z\psi_0+\frac{1}{2}z^2\left(\pa_z\psi_0\right)^2-\frac{1}{2z^2}(1+\sin\tau)\left(1-\frac{i\alpha}{z}\pa_z\psi_0\right)-2\eta i \left(z\pa_z\psi_0\right)^3=0.
\]
\end{remark}

\begin{figure}[h]
\begin{center}
\psfrag{t}{$\arctan\tet$}\psfrag{k1}{$i\kk$}\psfrag{k2}{$-i\kk$}\psfrag{D1}{
$\DD_{\kk,\tet}^{+,u}$}\psfrag{D2}{$\DD_{\kk,\tet}^{+,s}$}
\includegraphics[height=6cm]{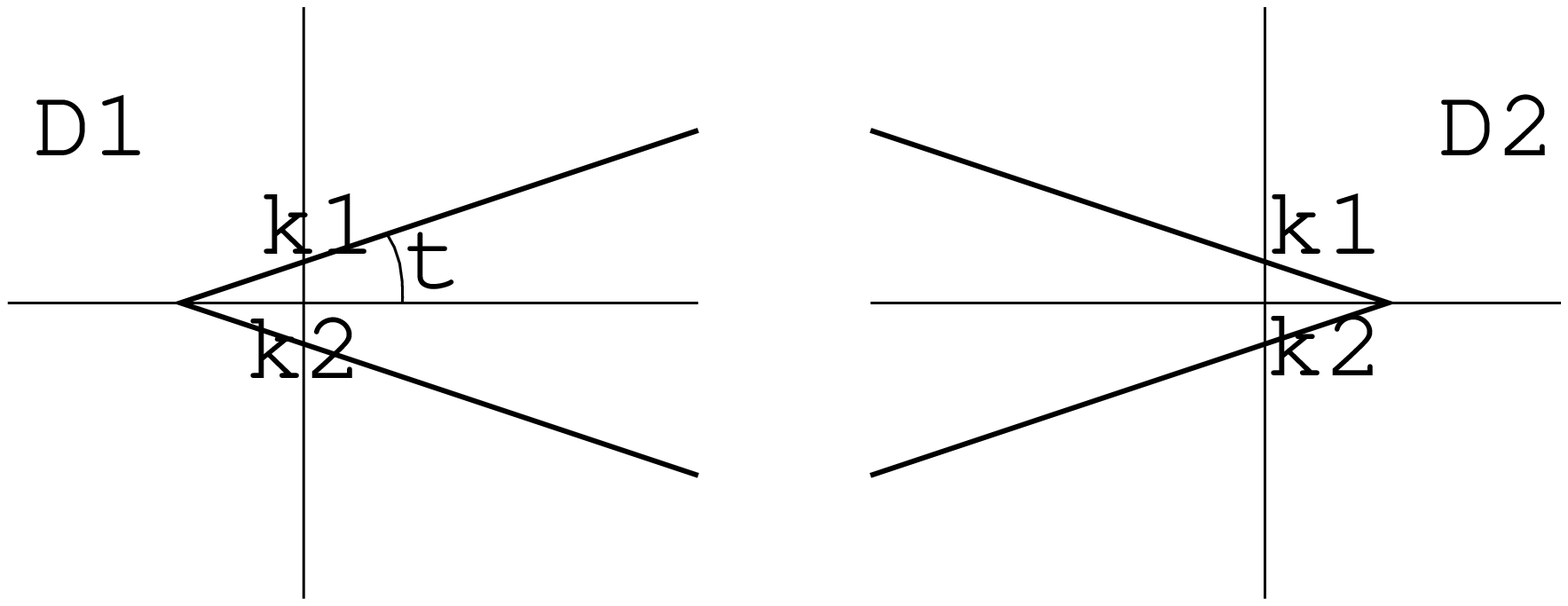}
\end{center}
\caption{\figlabel{fig:DomEqInner} The domains
$\DD_{\kk,\tet}^{+,u}$ and $\DD_{\kk,\tet}^{+,s}$ defined in
\eqref{def:Dominis:eqInner}.}
\end{figure}

In \cite{Baldoma06}, the author studies an inner equation slightly different  from \eqref{eq:HJEqInner} (see also \cite{OliveSS03}). The main differences between the inner equation considered in  \cite{Baldoma06}  and \eqref{eq:HJEqInner} is that in  \cite{Baldoma06}, the author assumes $\wt G=0$ and assumes also polynomial dependence on $\pa_z\psi_0$ whereas \eqref{eq:HJEqInner} can have the full jet. Nevertheless, one can easily see that using the techniques considered in \cite{Baldoma06} one can obtain exactly the same results obtained in that paper for inner equations of the form \eqref{eq:HJEqInner} since $\wt G$ is cubic and therefore, for $z$ big enough, is smaller than the other terms. We state these results in the next theorem, but first we introduce certain domains. The solutions of \eqref{eq:HJEqInner} are studied in the complex domains
\begin{equation}\label{def:Dominis:eqInner}
\begin{split}
\DD_{\kk,\tet}^{+,u}&=\left\{z\in\CC; \left|\Im z\right|>\tet \, \Re
z+\kk\right\}\\
\DD_{\kk,\tet}^{+,s}&=\left\{z\in\CC; -z\in \DD_{\kk,\tet}^{+,u}\right\}
\end{split}
\end{equation}
for $\kk>0$ and $\tet>0$. Let us observe that, for any $\C>0$,
$\DD_{\kk,\C}^{\inn,+,\ast}\subset\DD_{\kk, \tan\beta_2}^{+,\ast}$
for $\ast=u,s$. Nevertheless, since through the proof we will have
to change the slope  of the domains $\DD_{\kk,\tet}^{+,\ast}$, we
start with a  slope $\tet=\tan \beta_2/2$. The difference between the stable and unstable manifolds of the
inner equation is studied in the intersection domain
\begin{equation}\label{def:DominisInterseccio:eqInner}
\RRR_{\kk,\tet}^{+}=\DD_{\kk,\tet}^{+,u}\cap \DD_{\kk,\tet}^{+,s}
\cap \left\{z\in \CC; \Im z<0\right\}.
\end{equation}
\begin{figure}[h]
\begin{center}
\psfrag{k}{$-i\kk$}\psfrag{D1}{$\DD_{\kk,\tet}^{+,u}$}\psfrag{D2}{$\DD_{\kk,\tet
}^{+,s}$}\psfrag{R}{$\RRR_{\kk,\tet}^{+}$}
\includegraphics[height=6cm]{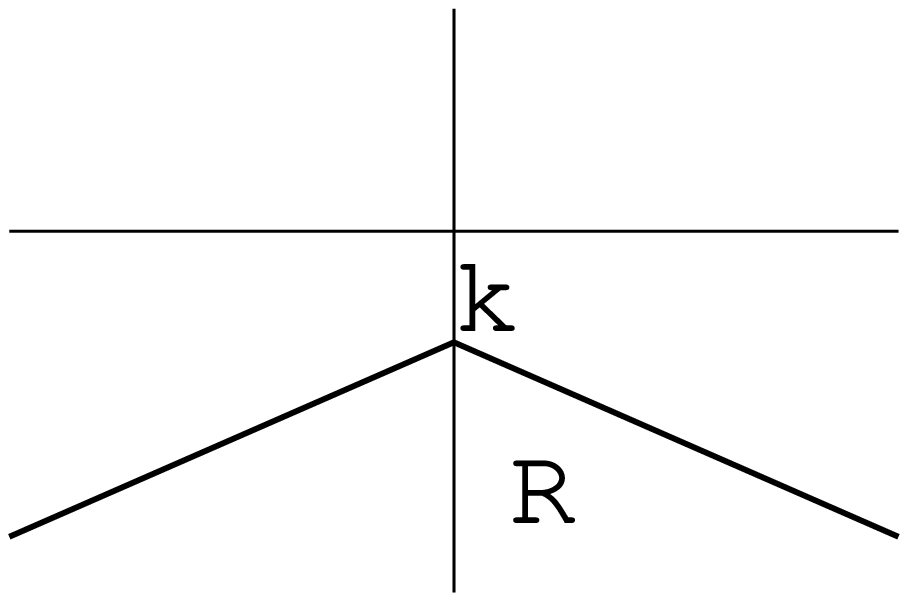}
\end{center}
\caption{\figlabel{fig:DomEqInnerInter} The domain
$\RRR_{\kk,\tet}^{+}$
 defined in \eqref{def:DominisInterseccio:eqInner}.}
\end{figure}

\begin{theorem}\label{th:InnerImma}
For any $\hmu\in B(\hmu_0)$
the following statements are satisfied:
\begin{enumerate}
\item There exists $\kk_2>0$ such that equation \eqref{eq:HJEqInner} has
solutions $\psi_0^{\ast}:
\DD_{\kk_2,\tan \beta_2/2}^{+,\ast}\times\TT_\sigma\rightarrow\CC$,
$\ast=u,s$, of the form
\begin{equation}\label{eq:SolucioInnerHJ}
\psi_0^{u,s}(z,\tau)=-\frac{1}{z}+\hmu\overline\psi_0^{u,s}
(z,\tau)+K^{u,s},\quad K^{u,s}\in\CC
\end{equation}
where $\overline\psi_0^{u,s}$ are analytic functions in all their
variables. Moreover, the derivatives of $\overline\psi_0^{u,s}$ are
uniquely determined by the condition
\[
\sup_{(z,\tau)\in
\DD_{\kk_2,\tan \beta_2/2}^{+,\ast}\times\TT_\sigma}\left|z^{3}\pa_z\ol\psi_0^{\ast
}(z,\tau)\right|<\infty
\]
for $\ast=u,s$.  In fact, one can choose $\ol\psi_0^{u,s}$ such that
\[
\sup_{(z,\tau)\in
\DD_{\kk_2,\tan \beta_2/2}^{+,\ast}\times\TT_\sigma}\left|z^{2}\ol\psi_0^{\ast}(z,
\tau)\right|<\infty
\]
for $\ast=u,s$.
\item There exists  $\kk_3>\kk_2$, analytic functions
$\left\{\chi^{[k]}(\hmu)\right\}_{k\in\ZZ^-}$
defined on $B(\hmu_0)$ and
$\zeta:\RRR_{\kk_3,\tan \beta_2}^{+}\times\TT_\sigma\rightarrow \CC$ such
that two solutions $\psi_0^{u,s}$ of equation \eqref{eq:HJEqInner}
of the form  given in \eqref{eq:SolucioInnerHJ} with $K^u=K^s$,
satisfy
\[
\left(\psi^u_0-\psi^s_0\right)(z,\tau)= \hmu\sum_{k<0}
\chi^{[k]}(\hmu)e^{ik\left(z-\tau+\hmu \zeta(z,\tau)\right)}.
\]
Moreover, the function $\zeta$  is of the form
\[
\zeta(z,\tau)=-\hmu b\ln z+\ol \zeta(z,\tau),
\]
where $b$ is the constant  defined in \eqref{def:Constantb} and
$\ol \zeta$ satisfies
\begin{equation*}
 \sup_{(z,\tau)\in
\RRR_{\kk_3,\tan \beta_2}^{+}\times\TT_\sigma}\left|z \ol
\zeta(z,\tau)\right|<\infty.
\end{equation*}
\end{enumerate}
\end{theorem}

The proof of this theorem follows the same lines as the theorem proved in \cite{Baldoma06}. Moreover, it can be easily seen that the
analytic functions  $\left\{\chi^{[k]}(\hmu)\right\}_{k\in\ZZ^-}$ are entire. 


To have a better knowledge of the parameterizations of the invariant
manifolds in the inner domains $\DD_{\kk, \C}^{\inn, +,\ast}$,
$\ast=u,s$ in \eqref{def:DominisInnerEnz}, we need to compare the
parameterizations $\psi^{u,s}$, which are solutions of
\eqref{eq:HJGeneralInner} with $\psi_0^{u,s}$ which are solutions of
\eqref{eq:HJEqInner} and have been given in Theorem
\ref{th:InnerImma}. Recall that $\DD_{\kk,\C}^{\inn,+,\ast}\subset\DD_{\kk, \tan\beta_2}^{+,\ast}$, 
 $\ast=u,s$, and therefore $\psi_0^{u,s}$ are defined in these domains.

We state the next theorem for the unstable invariant manifold. The
stable manifold satisfies analogous properties.

\begin{theorem}\label{th:MatchingHJ}
Let $\ga\in (0,1)$, the constants $\kk_1$ and  $\kk_3$  defined in Theorems
\ref{th:ExtensioFinal} and \ref{th:InnerImma},
$\C_1>0$ and  $\eps_0>0$
 small enough and $\kk_4>\max\{\kk_1,\kk_3\}$ big enough, which might depend on
the previous constants.
Then, for $\eps\in (0,\eps_0)$ and $\hmu\in B(\hmu_0)$, there exists a
constant $b_{2}>0$ such that for $(z,\tau)\in \DD_{\kappa_4,
\C_1}^{\inn,+,u}\times\TT_\sigma$,
\[
\left|\pa_z\psi^u(z,\tau)-\pa_z\psi_0^u(z,\tau)\right|\leq
\frac{b_{2}\eps^{\frac{2}{M}}}{\left|z\right|^{2-\frac{2}{M}}},
\]
where $\ga$ enters in the definition of
$\DD_{\kappa_4, \C_1}^{\inn,+,u}$, $\psi^u_0$ is given in Theorem
\ref{th:InnerImma}, $\psi^{u}$ is the scaling of the generating
function $T^u$ given in \eqref{eq:FuncioGeneradoraInner} and $M$ has been given in \textbf{HP2}.
\end{theorem}
The proof of this theorem follows the same lines as the proof of Theorem 4.16 of \cite{BaldomaFGS11}.

\subsection{Study of the difference between the invariant
manifolds}\label{sec:Difference}

Once we have obtained parameterizations of the invariant manifolds
of the form \eqref{eq:ParameterizationHJ} in the domains
$D^s_{\kk_1,d_1}$ and $D^u_{\kk_1,d_1}$ and studied their first order
approximation near the singularities, the next step is to study their
difference. To this end, we define
\begin{equation*}
\Delta(u,\tau)=T^u(u,\tau)-T^s(u,\tau)
\end{equation*}
in the domain 
$R_{\kk,d}=D^s_{\kk,d}\cap D^u_{\kk,d}$
which is defined as
\begin{equation}\label{def:DominiRaro:Interseccio}
\begin{split}
R_{\kk,d}=\left\{u\in\CC; \right.&\left. |\Im u|<\tan\beta_2\Re
u+a-\kk\eps, |\Im u|>\tan\beta_2\Re u+a-d , \right. \\
&\left.|\Im u|<-\tan \beta_1\Re u+a-\kk\eps\right\},
\end{split}
\end{equation}
and can be seen in Figure \ref{fig:BoomInter}.

We recall that $p_0(u)\neq 0$ if $u\in R_{\kk,d}$ and hence  in this domain we can
use the functions $T^{s,u}$.
\begin{figure}[H]
\begin{center}
\psfrag{D6}{$\wt D^{\out,s}_{\rr,d,\kk}$} \psfrag{D5}{$\wt
D^{\out,u}_{\rr,d,\kk}$}
\psfrag{b1}{$\beta_1$}\psfrag{b2}{$\beta_2$}
\psfrag{a1}{$ia$}\psfrag{a2}{$-ia$}\psfrag{a3}{$i(a-d)$}\psfrag{a4}{
$i(a-\kk\eps)$}
\psfrag{D3}{$R_{\kk,d}$}\psfrag{D}{$D^{s}_{\kk,d}$}\psfrag{D1}{$D^{u}_{\kk,d}$}
\includegraphics[height=6cm]{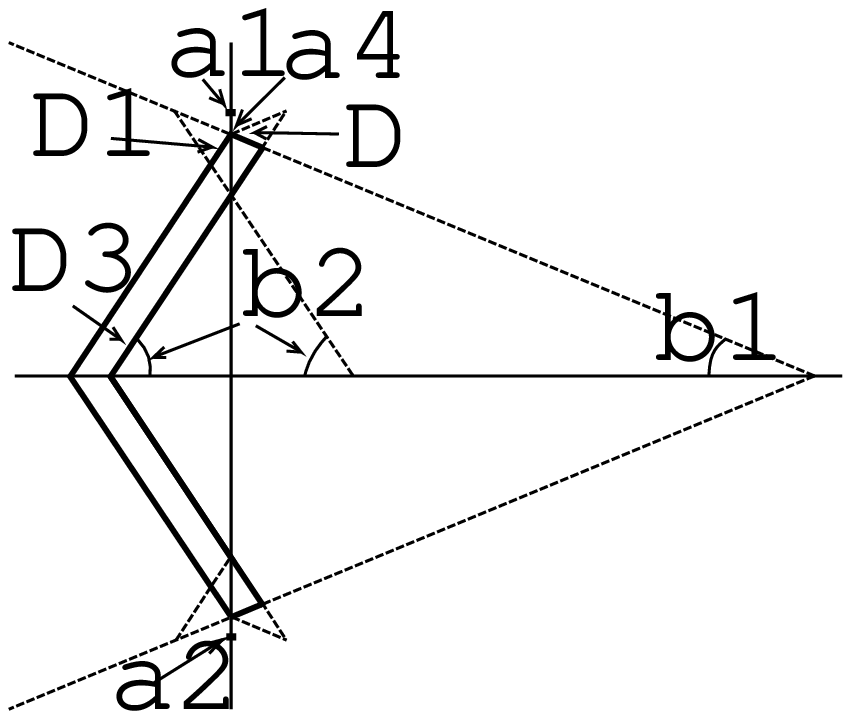}
\end{center}
\caption{\figlabel{fig:BoomInter} The domain $R_{\kk,d}$  defined in
\eqref{def:DominiRaro:Interseccio}.}
\end{figure}

Subtracting equation \eqref{eq:HamJacGeneral} for $T^u$ and $T^s$, one can easily see that $\Delta$ satisfies the partial differential equation
\begin{equation}\label{eq:DifferencePDE}
\widetilde\LL_\eps\xi=0,
\end{equation}
where 
\begin{equation*}
\widetilde\LL_\eps=\eps^{-1}\pa_\tau+(1+A(u,\tau))\pa_u
\end{equation*}
with 
\[
\begin{split}
A(u,\tau)=&\dps\frac{1}{2p_0^2(u)}\left(\pa_u T_1^u(u,\tau)+\pa_u
T_1^s(u,\tau)\right)\\
&+\frac{\hmu}{p_0(u)}\int_0^1 \pa_p \wh H_1
\left(q_0(u), p_0(u) +\frac{s\pa_u T_1^u(u,\tau)+(1-s)\pa_u
T_1^s(u,\tau)}{p_0(u)}, \tau\right)\, ds,
\end{split}
\]
where $\wh H_1$ is the function defined in
\eqref{def:ham:ShiftedOP:perturb} and
$T^{u,s}(u,\tau)=T_0(u)+T_1^{u,s}(u,\tau)$ with $\pa_u
T_0(u)=p_0^2(u)$ and $T_1^{u,s}$ are given in Theorem  \ref{th:ExtensioFinal}.

Following \cite{Baldoma06} (see also \cite{BaldomaFGS11}), to obtain the asymptotic expression of
the difference $\Delta$, we take advantage from the fact that it is
a solution of the homogeneous linear partial differential equation
\eqref{eq:DifferencePDE}. In \cite{Baldoma06} it is seen that if
\eqref{eq:DifferencePDE} has a solution $\xi_0$ such that
$(\xi_0(u,\tau),\tau)$ is injective in $R_{\kk,d}\times\TT_\sigma$,
then any solution of equation \eqref{eq:DifferencePDE} defined in
$R_{\kk,d}\times\TT_\sigma$ can be written as
$\xi=\Upsilon\circ\xi_0$ for some function $\Upsilon$.

Following this approach, we begin by looking  for a solution of the form
\[
\xi_0(u,\tau)=\eps\ii u-\tau+\CCC(u,\tau)
\]
where $\CCC$ is a function $2\pi$-periodic in $\tau$, such that
$(\xi_0(u,\tau),\tau)$ is injective in $R_{\kk,d}\times\TT_\sigma$. Moreover, to study the look for the first order of $\Delta$, we need to compare such function $\CCC$  with the function $\zeta$ obtained in Theorem \ref{th:InnerImma} in the inner domains \eqref{def:DominisInnerEnu}. 


\begin{theorem}\label{th:CanviFinal}
Let us consider  the constants $d_1>0$ defined in
Theorem \ref{th:ExtensioFinal} and $\kk_4>0$ in Theorem \ref{th:MatchingHJ}, $\eps_0>0$ small enough and $\kk_5>\kk_4$ big enough, which
might depend on the previous constants. Then, for $\eps\in
(0,\eps_0)$, $\mu\in B(\mu_0)$ and any $\kk\geq\kk_5$ such that
$\eps\kk<a$, there exists a real-analytic function
$\CCC(u,\tau):R_{\kk,d_1}\times\TT_\sigma\rightarrow\CC$ such that
\begin{itemize}
 \item $\xi_0(u,\tau)=\eps\ii u-\tau+\CCC(u,\tau)$ is solution of
\eqref{eq:DifferencePDE} and
\[
\left(\xi_0(u,\tau),\tau\right)=\left(\eps\ii
u-\tau+\CCC(u,\tau),\tau\right)
\]
is injective.
\item There exists a constant $b_{3}>0$ independent of $\mu$,
$\eps$ and $\kk$, such that for $(u,\tau)\in
R_{\kk,d_1}\times\TT_\sigma$,
\[
\begin{split}
\left|\CCC(u,\tau)\right|&\leq b_{3}\left|\hmu\right|\ln\left|u^2+a^2\right|\\
\left|\pa_u\CCC(u,\tau)\right|&\leq
\frac{b_{3}\left|\hmu\right|}{\left|u^2+a^2\right|}.
\end{split}
\]
\item Moreover, if we consider the constant $\C_1$  given in Theorem
\ref{th:MatchingHJ} and $\ga$ satisfying
\[
\frac{M}{2+M}<\ga<1,
\]
where $M$ has been defined in Hypothesis \textbf{HP2}. Then, there exist a constant $C(\hmu)$
defined for $\hmu \in B(\hmu_0)$ depending
real-analytically in
$\hmu$ and a constant $b_{4}>0$ such that $|C(\hmu)|\leq
b_{4}|\hmu|$ and, if $(u,\tau)\in \left(
D_{\kk,\C_1}^{\inn,+,u}\cap
D_{\kk,\C_1}^{\inn,+,s}\right)\times\TT_\sigma$ for any $\kk>\kk_5$,
\[
\left|\CCC(u,\tau)-C(\hmu)+\hmu b \ln\eps -\hmu
\zeta\left(\eps\ii(u-ia),\tau\right)\right|\leq
\frac{b_{4}|\hmu|\eps}{|u-ia|},
\]
where  $\zeta$ is the function given in Theorem
\ref{th:InnerImma} and $b$ is the constant defined in
\eqref{def:Constantb} respectively.
\end{itemize}
\end{theorem}
The proof of this theorem follows the same lines as the proofs of Theorem 4.21  and Proposition 4.22 of \cite{BaldomaFGS11}.

As we have explained, since $\Delta=T^u-T^s$ is a solution of the
same homogeneous partial differential equation as $\xi_0$ given in
Theorem \ref{th:CanviFinal}, there exists a function
$\Upsilon$ such that $\Delta=\Upsilon\circ\xi_0$, which gives
\begin{equation*}
\Delta(u,\tau)=\Upsilon\left(\eps\ii u-\tau+\CCC(u,\tau)\right).
\end{equation*}
Since $\Delta$ is $2\pi$-periodic in $\tau$,
 the function $\Upsilon$ is also $2\pi$-periodic in its variable.
Therefore, considering the Fourier series of $\Upsilon$ we obtain
\begin{equation}\label{def:Delta:Fourier}
\Delta(u,\tau)=\sum_{k\in\ZZ}\Upsilon^{[k]} e^{ik\left(\eps\ii
u-\tau+\CCC(u,\tau)\right)}.
\end{equation}

Now we  find the first asymptotic term of $\Delta$ which
is strongly related with
$(\psi^{u}_0-\psi^{s}_{0})(\eps^{-1}(u-ia),\tau)$, where
$\psi_0^{u,s}$ are the solutions of the inner equation given in Theorem
\ref{th:InnerImma}. We introduce the auxiliary function
\begin{equation*}
\Delta_0^+(u,\tau)=\sum_{k<0}\Upsilon_0^{[k]}e^{ik\left(\eps\ii
u-\tau+\CCC(u,\tau)\right)}
\end{equation*}
with
\begin{equation}\label{def:Dif0:coefsFourier}
\Upsilon_0^{[k]}=\frac{C_+^2\hmu}{\eps}\chi^{[k]}(\hmu)
e^{-{\dps\tfrac{|k|a}{\eps}}-i|k|(-C(\hmu)+\hmu b\ln\eps)}
\end{equation}
where $\left\{\chi^{k}(\hmu)\right\}_{k<0}$ are the coefficients
given in Theorem \ref{th:InnerImma} and  $C(\hmu)$ and $b$ are
the constants obtained in Theorem \ref{th:CanviFinal} and Theorem \ref{th:InnerImma} respectively. The scaling $C_+^2/\eps$
comes from the inner change in \eqref{eq:FuncioGeneradoraInner}.

We also introduce 
\begin{equation*} \Delta_0^-(u,\tau) =
\sum_{k>0}{\Upsilon_0^{[k]}}e^{ik\left(\eps\ii
u-\tau+\CCC(u,\tau)\right)}
\end{equation*}
with
\begin{equation}\label{def:Dif0:coefsFourier:positiu}
\Upsilon_0^{[k]}=\frac{\ol C_+^2\hmu}{\eps }\ol
\chi^{[-k]}(\hmu) e^{-{\dps\tfrac{|k|a}{\eps}}+i|k|(-\ol
C(\hmu,\eps)+\hmu \ol b\ln\eps)}.
\end{equation}
The function $\Delta_0^-(u,\tau)$  corresponds to the difference of
the solutions of the inner equation close to $u=-ia$. Taking $\tau, \hmu\in \RR$, $\Delta_0^{-}$
is nothing but the complex conjugate of $\Delta_0^{+}$. In fact, as
we know that $\Delta$ is a real analytic function in the $u$
variable for real values of $\hmu,\tau$, we can define $\Delta_0^-$
as the function that satisfies that $\Delta_0=\Delta_0^+ +
\Delta_0^-$ is also a real analytic function $\tau, \hmu\in \RR$. 
We will see that the first order of $\Delta$ is given by
\begin{equation*}
\Delta_0(u,\tau)=\Delta_0^+(u,\tau)+\Delta_0^-(u,\tau),
\end{equation*}
which can be written as
\begin{equation*}
\Delta_0(u,\tau)=\sum_{k\in\ZZ\setminus\{0\}}\Upsilon_0^{[k]}e^{ik\left(\eps\ii
u-\tau+\CCC(u,\tau)\right)},
\end{equation*}
where $\Upsilon_0^{[k]}$ are defined either by
 \eqref{def:Dif0:coefsFourier} and
\eqref{def:Dif0:coefsFourier:positiu}. For convenience we introduce $\Upsilon_0^{[0]}=0$. From
now on, we consider real values of
$\tau\in\TT_\sigma\cap\RR$.

\begin{theorem}\label{th:CotaDiferencia}
Let us consider  the mean value of $\Upsilon$, $\Upsilon^{[0]}$,
defined in \eqref{def:Delta:Fourier}, $s<2/M$, and $\eps_0>0$ small enough.
Then, there exists a constant $b_{5}>0$ such that for
$\eps\in(0,\eps_0)$ and $\hmu\in B(\hmu_0)\cap\RR$ and $(u,\tau)\in
\left(R_{s\ln(1/\eps),d_1}\cap\RR\right)\times\TT$, the following
statements are satisfied.
\[
\begin{split}
\left|\Delta(u,\tau)-\Upsilon^{[0]}- \Delta_0(u,\tau)\right|&\leq
\frac{b_{5}|\hmu|}{\eps|\ln\eps|}e^{-{\dps\tfrac{a}{\eps}}+\hmu \Im
b\ln\eps}\\
\left|\pa_u\Delta(u,\tau)-\pa_u \Delta_0(u,\tau)\right|&\leq
\frac{b_{5}|\hmu|}{\eps^{2}|\ln\eps|}e^{-{\dps\tfrac{a}{\eps}}+\hmu \Im
b\ln\eps}\\
\left|\pa_u^2\Delta(u,\tau)-\pa_u^2 \Delta_0(u,\tau)\right|&\leq
\frac{b_{5}|\hmu|}{\eps^{3}|\ln\eps|}e^{-{\dps\tfrac{a}{\eps}}+\hmu \Im
b\ln\eps}.\\
\end{split}
\]
\end{theorem}
The proof of this Theorem follows the same lines as the proof of Theorem 4.23 in \cite{BaldomaFGS11}. Note that, following the notation of \cite{BaldomaFGS11},  the case considered in Theorem \ref{th:CotaDiferencia} corresponds to  $r=1$ and $\ell=2$ in  \cite{BaldomaFGS11} and thus $\ell-2r=0$.

We observe that $\partial_u \Delta_0$ gives the correct asymptotic
prediction of $\partial_u \Delta$ if $\Upsilon^{[-1]}_0\neq 0$. In
fact, we only need this coefficient to give a simpler leading term
of the asymptotic formula. For this purpose let us define the
function
\begin{equation*}
f\left(\hmu\right)=C_+^2 \chi^{[-1]}\left(\hmu\right),
\end{equation*}
where $C_+$ is the constant defined in 
\eqref{eq:SepartriuAlPolTrig} and $\chi^{[-1]}(\hmu)$ is the
constant given in Theorem \ref{th:InnerImma}.  This function $f(\mu)$ is the one appearing in  the asymptotic formula for the area of the lobes \eqref{def:formulaArea:singular} given in Theorem \ref{th:Main}. From its definition just given, we can see that this constant $f(\mu)$  is essentially given by $\chi^{[-1]}(\hmu)$, that is, by the first harmonic Fourier coefficient of the difference between the solutions $\psi_0^{u,s}$ of the inner equation \eqref{eq:HJEqInner} (see Theorem \ref{th:InnerImma}). In particular, the zeros of $f(\hmu)$ correspond to the zeros of $\chi^{[-1]}(\hmu)$.

We
define
\begin{equation*}
\Delta_{00}(u,\tau)=\frac{2\hmu}{\eps^{2r-1}}e^{{\dps
-\tfrac{a}{\eps}}}\Re\left(f(\hmu)e^{-i\left(\hmu^2b\ln\eps-C(\hmu)\right)}
e^{-i\left({\dps
\tfrac{u}{\eps}}-\tau+\CCC(u,\tau)\right)}\right)
\end{equation*}
where the constant $b$  has been defined in \eqref{def:Constantb} and the function $\CCC$   and the constant 
$C(\hmu)$ have been defined in Theorem \ref{th:CanviFinal}.

\begin{corollary}\label{coro:Diferencia:Simple}
There exists a constant $b_{6}>0$ such that for
$\eps\in(0,\eps_0)$, $\hmu\in B(\hmu_0)\cap\RR$ and $(u,\tau)\in
\left(R_{s\ln(1/\eps),d_1}\cap\RR\right)\times\TT$, the following
statements are satisfied,
\[
\begin{split}
\left|\Delta(u,\tau)-\Upsilon^{[0]}- \Delta_{00}(u,\tau)\right|&\leq
\frac{b_{6}\left|\hmu\right|}{\eps|\ln\eps|}e^{-{\dps\tfrac{a}{\eps}}
+\hmu \Im b\ln\eps}\\
\left|\pa_u\Delta(u,\tau)-\pa_u \Delta_{00}(u,\tau)\right|&\leq
\frac{b_{6}\left|\hmu\right|}{\eps^{2}|\ln\eps|}e^{-{\dps\tfrac{a}{\eps}}
+\hmu \Im b\ln\eps}\\
\left|\pa_u^2\Delta(u,\tau)-\pa_u^2 \Delta_{00}(u,\tau)\right|&\leq
\frac{b_{6}\left|\hmu\right|}{\eps^{3}|\ln\eps|}e^{-{\dps\tfrac{a}{\eps}}
+\hmu \Im b\ln\eps}.\\
\end{split}
\]
\end{corollary}

To finish the proof of Theorem \ref{th:Main}, it is enough to proceed as in Section 4.8 of \cite{BaldomaFGS11}.

\section*{Acknowledgements}
M. G. acknowledges useful discussions with Inmaculada Baldom\'a, Yuri Fedorov, Ernest Fontich and Tere M. Seara. He  has been partially supported by the Spanish MCyT/FEDER grant MTM2009-06973 and the Catalan SGR grant 2009SGR859.  Part of this work was done while he was doing stays in  the University of Maryland at College Park, the Pennsylvania State University, the Fields Institute and Universitat Polit\`ecnica de Catalunya.  He wants to thank these institutions for their hospitality and support.

\bibliography{references}
\bibliographystyle{alpha}
\end{document}